\documentclass[11pt]{amsart}

\usepackage{a4wide}
\usepackage[T1]{fontenc}
\usepackage[latin1]{inputenc}
\usepackage{amssymb}
\usepackage{amsfonts}
\usepackage{latexsym}
\usepackage[dvips]{graphics}
\usepackage{psfrag}
\makeatletter

\newcommand{\cal}{\mathcal}
\newcommand{\1}{{\bold 1}}
\newcommand{\<}{\langle}
\newcommand{\R}{\mathbb{R}}
\newcommand{\C}{\mathbb{C}}

\newcommand{\MM}{\cal {M}}

\newcommand{\EE}{\cal {E}}
\newcommand{\HH}{\mathcal {H}}
\newcommand{\DD}{{\cal {D}}}
\newcommand{\N}{\mathbb{N}}
\newcommand{\te}{\theta}

\newcommand{\al}{\alpha}
\newcommand{\be}{\beta}
\newcommand{\si}{\sigma}
\newcommand{\la}{\lambda}
\newcommand{\ep}{\varepsilon}
\newcommand{\Om}{\Omega}
\newcommand{\om}{\omega}
\newcommand{\p}{\partial}
\newcommand{\ti}{\tilde}
\newcommand{\Ti}{\widetilde}

\newcommand{\To}{\longrightarrow}

\newcommand{\De}{\Delta}
\newcommand{\de}{\delta}
\newcommand{\ga}{\gamma}

\newcommand{\vide}{\emptyset}
\newcommand{\vphi}{\varphi}
\newcommand{\tr}{\mbox{tr}}
\newcommand{\re}{\mbox{Re}}

\newcommand{\Res}{\mbox{Res}}
\newcommand{\res}{\emph{Res}}

\def \Op{Op_{h}^{\om}}
\def \.{{\bf{\cdot}}}
\def \Proof{{\bf{Proof.}  }}
\newtheorem{prop}{Proposition}
\newtheorem{lem}{Lemma}
\newtheorem{thm}{Theorem}
\newtheorem{rem}{Remark}
\newtheorem{corl}{Corollary}
\newtheorem{defi}{Definition}

\begin{document}
\makeatother
\title[Dirac Resonances and Spectral Shift Function]{Resonances and Spectral Shift Function for the semi-classical Dirac operator\vspace{5mm}}
\author{Abdallah  Khochman}
\date{\today}
\email{Abdallah.Khochman@math.u-bordeaux1.fr}
\address{Math\'ematiques Appliqu\'ees, Universit\'e Bordeaux I, 351, cours de la Lib\'eration, 33405 Talence, France}
\maketitle


\begin{abstract}
We consider the self-adjoint operator $H=H_0+V$, where $H_0$ is the free semi-classical Dirac operator on $\R^3$. We suppose that the smooth matrix-valued potential $V=O(\<x\rangle^{-\de}),\;\de>0,$ has an analytic continuation in a complex sector 
 outside a compact. We define the resonances as the eigenvalues of the non-selfadjoint operator obtained from the Dirac operator $H$ by a complex distortions of  $\R^{3}$.
We establish an upper bound $O(h^{-3})$ for the number of resonances in any compact domain. For $\de>3$, a representation of the derivative of the spectral shift function $\xi(\la,h)$ related to the semi-classical resonances of $H$ and a local trace formula are obtained. In particular, 
if $V$ is an electro-magnetic potential, 
 we deduce a Weyl-type asymptotic of the spectral shift function. 
 As a by-product, we obtain an upper bound $O(h^{-2})$ for the number of resonances close to non-critical energy levels in domains of width $h$ and a Breit-Wigner approximation formula for the derivative of the spectral shift function.
\end{abstract}
\maketitle{{\it Keywords:} Semi-classical Dirac operator - Resonances - Trace formula - Spectral shift 
function - Weyl-type asymptotics - Breit-Wigner approximation.

{\it Mathematics classification:} 35B34 - 35P05 - 34L40 - 81Q20 - 81Q10. }
\section{Introduction}
$\;\;\;\;$The resonance theory for 
the Schr\"odinger equation has been developed  following  several approaches. Among them we can mention
the analytic dilation (see \cite{JAJC}) or the analytic distortion (see \cite{WH}) and, in the semi-classical regime, that related to the work of Helffer-Sj\"ostrand \cite{BHJS}. In \cite{BHAM} Helffer-Martinez  showed that the different definitions give  the same resonances when one can simultaneously apply them to an operator. For the three dimensional Dirac operator,  Seba \cite{PS} defined the resonances 
  as complex eigenvalues of the operator obtain by a complex dilation. Applying the approach of Helffer-Sj\"ostrand \cite{BHJS}, Parisse \cite{BP} has studied the Dirac resonances in the semi-classical regime, with some scaling functions. The last two works deal with analytic perturbations near the real axis.

The concept of the spectral shift function has been introduced by Lifshits \cite{IML} in connection with the problems in quantum statistics and solid physics. Thereafter, a mathematical theory of the spectral shift function has been constructed by Krein \cite{MK}. Moreover, in \cite{MBMK} Birman-Krein found a connection between the scattering theory and the theory of the spectral shift function. 
A detailed presentation of the theory of the spectral shift function can be found in \cite{DY}.
For a survey concerning the spectral shift function (SSF) for Schr\"odinger and Dirac operator  or the asymptotic expansion of this function, 
we refer to Robert \cite{DR2} and  to the references given there.

A representation of the derivative of the scattering phase in terms of the resonances has been established for the Schr\"odinger operators. Such representations have been successively  obtained by Melrose \cite{RM} for obstacle problems in the high energy case, by Petkov-Zworski \cite{VPMZ1}, \cite{VPMZ2} for "black box" scattering with compact perturbations in the classical and the semi-classical cases and by Bruneau-Petkov \cite{VBVP} 
for long-range perturbations in the semi-classical "black box" framework. The results in \cite{VBVP} have been generalized by Dimassi-Petkov \cite{MDVP} for non-semi-bounded Schr\"odinger type operators.
As a by-product, they prove a Weyl type asymptotic for the the scattering phase.
Moreover, Weyl asymptotic can also be obtained by representation of the derivative of the spectral shift function involving the trace of the cut-off resolvent (see Robert \cite{DR3}, Bruneau-Petkov \cite{VBVP1} and Nakamura \cite{SN}).

Concerning the Breit-Wigner approximation for the derivative of the spectral shift function in the Schr\"odinger case, same kind of results have been obtained in a particular semi-classical set-up by C.G\'erard-Martinez-Robert \cite{CGAMDR} for short range potentials on $\R^n$ and by Petkov-Zworski \cite{VPMZ2} for a general compactly supported perturbation (see also \cite{VBVP1}, \cite{JBJS}).\\

For Dirac operators, Bruneau-Robert \cite{VBDR} established an asymptotic expansion of the scattering phase $s(\la)$ and their derivatives in the high energy regime and in the semi-classical regime for $\la$ in a non-trapping energy interval. For an interval $I \subset ]-mc^2,mc^2[$ with non critical extremities, 
  Helffer and Robert in \cite{BHDR} gave an asymptotics of the number of the eigenvalues in $I$ for scalar potentials. Nevertheless, we are not aware of works which deal with the link between the derivative of the SSF and resonances for the semi-classical Dirac operators (in the spirit of Petkov-Zworski \cite{VPMZ2} and Bruneau-Petkov \cite{VBVP}), 
neither of papers giving the Weyl asymptotic of the spectral shift function for Dirac operators in any interval $I$. \\ 

The purpose of this work is to extend the definition of resonance for analytic perturbations outside a compact set. We define the resonances for the semi-classical Dirac operator as the discrete eigenvalues of the non-selfadjoint operator obtained from the Dirac operator $H$ by a general class of complex distortions of  $\R^{3}$. We prove that the resonances are independent of the  distortion (see Section \ref{secdefinition}). We establish an upper bound  
for the number of resonances in a compact domain $\Om$ (see Section \ref{majoration}). The second goal of this work is to obtain a meromorphic continuation of the derivative of the spectral shift function $\xi(\la,h)$ related to the resonances for the semi-classical Dirac operator (see Section \ref{secrep}). This last is closely related to trace formulae (see \cite{VBVP}, \cite{JS1}, \cite{JS2}, \cite{VPMZ1}, \cite{VPMZ2}, \cite{JSMZ}) and to resonance expansions (see \cite{STMZ}, \cite{NBMZ}). Thereafter, in the case where the potential is an electro-magnetic potential, we deduce a Weyl-type asymptotic of the spectral shift function (see Section \ref{secweyl}). As a by-product, we obtain  an upper bound $O(h^{-2})$ for the number of resonances close to non-critical energy levels in domains of width $h$ (see Subsection \ref{majorationh}), as well as a Breit-Wigner approximation for the derivative of the SSF (see Subsection \ref{secbreit}).
 \section{statement of the results}
 We consider the selfadjoint Dirac operator 
\begin{eqnarray}\label{eqH_0}H_0=-ich\sum_{j=1}^{3}\alpha_j \frac{\p }{\p x_j} +
\beta mc^2,\end{eqnarray}     
with domain $D(H_0)=H^1( \R^3)\otimes { \C}^4\subset {\cal H}
=L^2(\R^3)\otimes \C^4,$ where $h\searrow0$ is the semi-classical parameter, $m>0$ is the mass of the Dirac particle and $c$ is the 
speed of light.
 The quantities $\alpha_1,\ \alpha_2,\ \alpha_3$ and $\beta$ are $4\times4$ Dirac matrices satisfying the anti-commutation relations
\begin{eqnarray}\label{eqanti}\left\{\begin{array}{ll}\alpha_i\alpha_j+\alpha_j\alpha_i=2\delta_{ij}I_4,\,\,\,\,\,\,\,\,\,&\mbox{for}\,\,\,\,\,\,\,\,\,i,j=1,2,3,\\
\alpha_i\beta+\beta\alpha_i=0,\,\,\,\,\,\,\,\,\,&\mbox{for}\,\,\,\,\,\,\,\,\,i=1,2,3,
\end{array}\right.
\end{eqnarray}
and $\be^2=I_4$. Here $I_n$ is the $n\times n$ identity matrix. For example, we choose the standard (or Dirac-Pauli) representation of these matrices \[\al_i=\left(\begin{array}{cc}
0&\si_i\\
\si_i&0
\end{array}\right),\;\;\;\be=\left(\begin{array}{cc}
I_2&0\\
0&-I_2
\end{array}\right),\]
where $(\si_j)_{1\leq j\leq3}$ are the $2\times2$ Pauli matrices:
\[\si_1=\left(\begin{array}{cc}
0&1\\
1&0
\end{array}\right),\;\;\;\si_2=\left(\begin{array}{cc}
0&-i\\
i&0
\end{array}\right),\;\;\;\si_3=\left(\begin{array}{cc}
1&0\\
0&-1
\end{array}\right).\] 
\begin{rem}
Most calculations with Dirac matrices can be done without referring to a particular representation (see Appendix 1.A \cite[Chap. 1]{BT}).
\end{rem}
Let 
$H_1\!=\!H\!:=\!H_0+V$,
where $V$ is the multiplication operator by a $4\!\times\!4$-matrix potential $V$. We suppose that $V\in C^{\infty}(\R^3)$  and satisfies the following assumption\\\\
${\bf(A_V)}:\;\mbox{}$ $V$ \emph{is Hermitian on $\R^3$ and has an analytic extension in the sector} 
\begin{eqnarray}\label{eqsector}
C_{\epsilon,0}:=\{z\in \C^3,|\mbox{Im} (z)|\leq\,\epsilon\,|\mbox{Re} (z)|,\;\;|\re(z)|>R_0\},\;   \;for\;0<\epsilon <1.  
\end{eqnarray} 
 $\hspace{13mm}\;Moreover,\;for\;x\in C_{\epsilon,0}\;it\; satisfies$ 
 \begin{eqnarray}\label{eqAV}
\;\;\;\;\;\;\;\; \|V(x)\|=O(\<x\rangle^{-\de}),\;\;\;\de>0,\;\; \<x\rangle=(1+|x|^2)^{\frac12}\;\mbox{}.
 \end{eqnarray}
 The free Dirac operator $H_0$ has only essential spectrum $\si_{ess}(H_0)=]-\infty,-mc^2]\cup[mc^2,+\infty[$.
Under the assumption  ${\bf(A_V)}$ the operator $H_1$ is a selfadjoint operator. Using Weyl theorem, we have  $\si_{ess}(H_1)=\si_{ess}(H_0).$

For $\te\in D_{\epsilon}:=\{\te\in\C,|\te|\leq r_{\epsilon}:=\frac{\epsilon}{\sqrt{1+\epsilon^2}}\},$ 
 we denote  $$H_{1,\te}=H_{\theta}:=U_{\theta}H_0U_{\theta}^{-1}+U_{\theta}VU_{\theta}^{-1}=H_{0,\te}+U_{\theta}VU_{\theta}^{-1},$$ where $U_{\theta}$ is the one-parameter family of distortions defined below (see Section $\ref{secdistortion}$).

For $\te_0$ fixed in $D_{\epsilon}^+:=D_{\epsilon}\cap\{\te\in\C,\;\mbox{Im} (\te)\geq0\}$, we define  
  \[\Gamma_{\te_0}:=\{\pm c\sqrt{\frac{\lambda}{(1+\theta_0 )^2} + m^2c^2}\in\C,\;\; \lambda\in [0,+\infty[\},\vspace{-4mm}\]
 and
 \[S_{\te_0}:=\{ z\in\bigcup_{
\substack{\te\in D_{\epsilon}^+}}\!\!\Gamma_{\te};\;\;\mbox{arg}(1+\te)<\mbox{arg}(1+\te_0),\;\;\frac{1}{|1+\te|}<\frac{1}{|1+\te_0|}\;\}.\vspace{-4mm}\]
\\\\
The square root $\sqrt{z}$ is defined such that for $z\in\C\backslash]-\infty,0]$, Re$(\sqrt{z})>0$.\\

For $\te\in D_{\epsilon}^+$, $ \mbox{arg}(1+\te_0)\leq\mbox{arg}(1+\te),\;\;\frac{1}{|1+\te_0|}\leq\frac{1}{|1+\te|}$, we prove that the spectrum of $H_{\te}$ is discrete in $S_{\te}$ and independent of $\te$ in $S_{\te_0}$. This justifies the following definition
\begin{defi}\label{defres}
The resonances of $H$ in  $S_{\te_0}\cup\R$ are the 
eigenvalues of $H_{\te_0}$.
The multiplicity of a resonance $z$ is the geometric multiplicity of $z$ considered as an eigenvalue of $H_{\te_0}$.
We will denote $\res(H)$ the set of resonances.
\end{defi}
%
The most important advantage of this definition is that the resonances can be computed by solving a non-selfadjoint eigenvalue problem. 
\begin{rem}
The resonances of $H$ in $\{z\in\C;\;\emph{Re}(z)\in]-mc^2,mc^2 [\}$ are the real eigenvalues of $H$.\\

\end{rem} \vspace{-1cm}

\begin{figure}[ltbh]
{\large
\psfrag{S}[lt][][1][0]{$S_{\te_0}$}
\psfrag{0}[l][][1][0]{}
\psfrag{q}[l][][1][0]{}
\psfrag{2}[l][][1][0]{}
\psfrag{2}[l][][1][0]{}
\psfrag{mc}[c][][1][0]{$mc^2$}
\psfrag{-mc}[c][][1][0]{$-mc^2$}
\psfrag{W}[c][][1][0]{$\Omega$}
\psfrag{Re\(z\)}[c][][1][0]{$\mbox{Re}(z)$}
\psfrag{Im\(z\)}[c][][1][0]{$\mbox{Im}(z)$}
\psfrag{G}[lt][][1][0]{$\Gamma_{\te_0}$}
\psfrag{Fig.1.}[c][][1][0]{Fig. The set $S_{\te_0}$}
\psfrag{Resonances}[l][][1][0]{Resonances}
\hspace{-1.7cm}\includegraphics{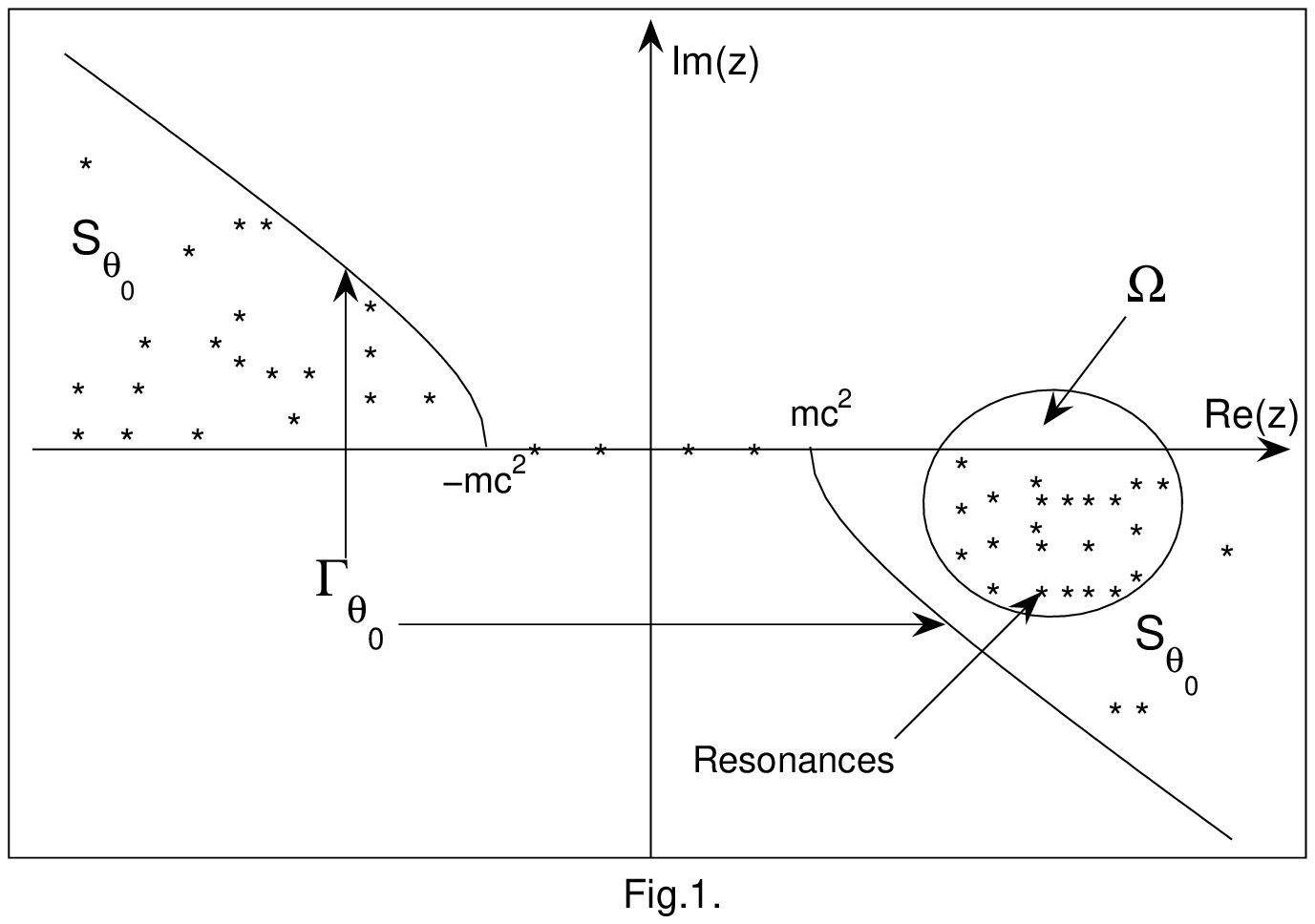}
}
\end{figure}


Now, we would like to find symmetry properties  which 
makes it  possible to limit our study of the resonances in a domain $\Om$  which satisfies ${\bf(A_{\Om}^+)}$, with\\ \\
${\bf(A_{\Om}^{\pm})}$: $\Om$ $is\;an\; open\;simply\;connected\;and\;relatively\; compact\;subset\;of\;\{z\in\C;\pm\mbox{Re} (z)>mc^2\}\;such \; that$ $\Om\cap\{\pm\mbox{Im} (z)>0\}\neq\emptyset\;\;$ $and$  $\;there\;exists\;\te_0\in D_{\epsilon}^+ 
  \;such\; that \;\; {\overline \Om}\cap\Gamma_{\te_0}=\emptyset $.  
 \begin{prop}\label{propimmaginaire}
Let $H^-$ be the selfadjoint Hamiltonian $$H^-=H_{0}-U_c{\overline {V}(x)}U_c^{-1},$$
where $U_c=i\be\al_2$ is an unitary matrix $4\times4$ and ${\overline V}$ is the conjugate of $V$. 
 Then the following assertions are equivalents:
 \begin{itemize} 
\item[\emph{(i)}] The complex value $z$  is a resonance of $H$.
\item[\emph{(ii)}] The symmetric of the conjugate $-{\bar z}$ is a complex eigenvalue of $U_{\bar\te}H^-U_{\bar\te}^{-1}$.
\end{itemize}
Moreover, the multiplicity of $z$ is equal to the multiplicity of $-{\bar z}$ considered as an eigenvalue of $U_{\bar\te}H^-U_{\bar\te}^{-1}$.
\end{prop}
\begin{prop}\label{propreelle}
Let ${\overline H}$ be the selfadjoint Hamiltonian $${\overline H}=-ic\sum_{j=1}^3\al_j'\p_{x_j}+\be'mc^2+{\overline V},$$ where $\al_1'=-\al_1,\;\al_2'=\al_2,\;\al_3'=-\al_3,\;\be'=\be$ 
 are matrices which satisfy the anti-commutation relations (\ref{eqanti}) and ${\overline V}$ is the conjugate of $V$. Then the following assertions are equivalents:
\begin{itemize} 
\item[\emph{(i)}] The complex value $z$  is a resonance of  $H$.  
\item[\emph{(ii)}] The conjugate $\bar z$ is a complex eigenvalue of $U_{\bar\te}{\overline H}U_{\bar\te}^{-1}$.
\end{itemize}
Moreover, the multiplicity of $z$ is equal to the multiplicity of ${\bar z}$ considered as an eigenvalue of $U_{\bar\te}{\overline H}U_{\bar\te}^{-1}$.
\end{prop}
 %
Using the same type of approach as in \cite{JS2}, we construct an operator ${\widehat H}_{j,\te}$ for $j=0,1$, so that 
 \[{\widehat H}_{j,\te}-{ H}_{j,\te}=K_j=O(1),\;\mbox{   has finite rank } O(h^{-3}),\]
 and $\|({\widehat H}_{j,\te}-z)^{-1}\|=O(1),$ uniformly for $z\in\,{\overline \Om}$ (see Subsection \ref{Appendix A}).\\

Using this construction we establish an upper bound  of the number of resonances: 
\begin{thm}{\bf(Upper bound)}\label{thmmajoration} Assume that $V$ satisfies the assumption ${\bf(A_V)}$ with $\de>0$. Let $\Om$ be a  complex domain satisfying the assumption ${\bf (A_{\Om}^{\pm})}$, then 
\[\#\res(H)\cap\Om \leq C(\Om)h^{-3}.\]
\end{thm}
For a pair of self-adjoint operators $(H_0,H_0+V)$ where $V$ satisfies the assumption $({\bf A_V})$ with $\de>3$, (see \cite{VBDR}, \cite{DR3}, \cite{DR2}), the spectral shift function $\xi(\la,h)$ is a distribution in ${\cal D}'(\R)$ such that its derivative satisfies: 
\begin{eqnarray}\label{defxi}\<\xi'(\la,h),f(\la)\rangle_{{\cal D}'(\R), {\cal D}(\R)}=\tr(f(H_1)-f(H_0)),\;\;\;f(\la)\in C_0^{\infty}(\R).\end{eqnarray}

By the Birman-Krein theory, the SSF is in $L_{loc}^1(\R)$ and coincides with the scattering phase (see \cite{DY}).
 
Our principal result  is a meromorphic continuation of the derivative of the spectral shift function $\xi(\la,h)$.
\begin{thm}\label{thmxi}{\bf(Representation formula)} Assume that $V$ satisfies the assumption ${\bf(A_V)}$ with $\de>3$. Let $\Om$ be a complex domain satisfying the assumption ${\bf(A_{\Om}^{\pm})}$ and $W\Subset \Om$ be an open simply connected set which is symmetric with respect to the real axis. Assume that 
$I=W\cap \R$ is an interval. Then for all $\la\in I$  we have the representation:
\begin{eqnarray}\label{eqderive}
\xi'(\la,h)=\frac{1}{\pi}\emph{Im} \,r(\la,h)+\sum_{\begin{array}{c}
\scriptstyle w\in \res(H_1)\cap\Om \\
\scriptstyle \emph{Im}\,  w\neq0
\end{array}}\frac{-\emph{Im}\,  w}{\pi|\la-w|^2}+\sum_{w\in \res(H_1)\cap I}\de_w(\la),
\end{eqnarray}
where $r(z,h)=g(z,h)-{\overline {g(\bar z,h)}}$, $g(z,h)$ is a holomorphic function in  $\Om$ which satisfies the following estimate: 
\begin{eqnarray}
|g(z,h)|\leq C(W)h^{-3}, \;\;z\in W,
\end{eqnarray}
with $C(W)>0$  independent of  $h\in]0,\;h_0]$. Here $\de_w(\cdot)$ is the Dirac 
mass at $w \in \R$. 
\end{thm}
\begin{rem} This theorem can be extended to the pair operators $(H_0+V_1$, $H_0+V_2)$, where $V_1,\;V_2$ are two $4\times4$ Hermitian potential matrices satisfying ${\bf(A_V)}$ with $\de>0$ and $V=V_2-V_1$ satisfies the assumption ${\bf(A_V)}$ with $\de>3$ (see Theorem \ref{thmxig}).
\end{rem}
As a corollary of the last theorem, we have a Sj\"ostrand type local trace formula (see Theorem \ref{thmtrace}).\\ 

 Now we discuss a Weyl type asymptotic of the spectral shift function $\xi(\la,h)$  in the case where $V$ is an electro-magnetic potential  
\begin{eqnarray}\label{eqV}
V(x)=e(-\al\cdot A+v)(x)=-\sum_{j=1}^3\al_j\cdot eA_j(x)+ev(x),\;\;
\end{eqnarray}
satisfying the assumption ${\bf(A_V)}$ with $\de>3$. Here $e<0$ is the charge of the Dirac particle. We assume that, the electric potential $v(x)=\left(
\begin{array}{cc}
v_+(x)I_2 & 0 \\
0 & v_-(x)I_2
\end{array}
\right)$ where $v_+,\ v_-$ are $C^\infty$ scalar  functions satisfying  \begin{eqnarray}\label{eqv+v-}
|e(v_+-v_-)(x)|<2mc^2,
\end{eqnarray}
and $A=(A_1,\ A_2,\ A_3)$ is a magnetic vector potential where $A_1,\ A_2,\ A_3$ are $C^\infty$ scalar  functions.

For any $(x,\xi)\in \R^6$,
 the semi-classical symbols of $H_\nu$, $\nu=0,1$ are the matrices
\begin{eqnarray}\label{eqDnu}
\DD_{\nu}(x,\xi)=\al\cdot(c\xi-\nu eA(x))+\be mc^2+\nu ev(x),\;\;\;\al=(\al_1,\al_2,\al_3),
\end{eqnarray}
which is Hermitian and has two eigenvalues
\begin{eqnarray}\label{eqH+-}
H_{\nu}^{\pm}(x,\xi)=\pm \left(\!|c\xi-\nu eA(x)|^2+\left(mc^2+\nu\frac{ e}{2}(v_+-v_-)\right)^2\!\right)^{\frac12}\!\!\!+\nu\frac{ e}{2}(v_++v_-),
\end{eqnarray}
of multiplicity two. The function $H_1^{+}(x,\xi)$ is the Hamiltonian for a relativistic classical particle and the other $H_1^{-}(x,\xi)$ can be considered as the one for its anti-particle (see \cite{PD}, \cite{BT}, \cite{KY}). Moreover, from (\ref{eqv+v-}) the two Hamiltonians $H_1^{\pm}(x,\xi)$ are smooth functions.

 For $\nu=0,1,$ the matrix
\begin{eqnarray}\label{eqP+-}
\Pi_\nu^{\pm}(x,\xi)=\frac12\left(1+\frac{\DD_\nu(x,\xi)-\nu ev(x)+\nu\frac{e}{2}(v_+-v_-)\be}{H_\nu^{\pm}(x,\xi)-\nu\frac{e}{2}(v_++v_-)}\right),
\end{eqnarray}
is the orthogonal projection on the eigenspace ${\EE}_\nu^{\pm}(x,\xi)$ of $\DD_\nu(x,\xi)$ corresponding to the eigenvalue $H_\nu^{\pm}(x,\xi).$
\begin{defi}\label{defnoncritical}
A real $\la$ is a noncritical energy level for $H_1$ if for all $(x,\xi)\in\R^6$, with $H_{1}^{\pm}(x,\xi)\!=\!\la,$ we have $\nabla_{x,\xi}H_{1}^{\pm}(x,\xi)\neq0.$
\end{defi}
\begin{thm}\label{thmweyl}{\bf(Weyl formula)}
Assume that the potential $V$ is an electro-magnetic potential given by $(\ref{eqV})$ and satisfying the assumption ${\bf(A_V)}$ with $\de>3$. For all $\la,\;\la_1$ noncritical energy levels for $H_1$ such that $\pm mc^2\not\in\; ]\la_1,\la[$ and $h\in ]0,h_0[$, 
we have the asymptotic expansion 
\begin{eqnarray}\label{weyl}
\xi(\la,h)-\xi(\la_1,h)=w(\la,\la_1)h^{-3}+O(h^{-2}).
\end{eqnarray}
Here the $O(h^{-2})$ is uniform for $\la$ (resp.$\la_1$) in a small interval $I_2$ (resp. $I_1$). The first term $w(\la,\la_1)\;\in C^\infty(I_2\times I_1)$ is given by
$$w(\la,\la_1)=w(\la)-w(\la_1),$$ with:
\begin{eqnarray}\label{eqpremierterm}
\;\;\;\;\;\;\;w(\la)=\frac{1}{3\pi^2}\int_{\R^3}\!W_+(\la,v_+,v_-)-W_+(\la,0,0)
- W_-(\la,v_+,v_-)+W_-(\la,0,0)\;\; dx
\end{eqnarray}
where   $W_\pm(\la,a,b)=\left(\!\left(\la-\frac{e(a+b)}{2}\right)_{\pm}^2-\left(mc^2+\frac{e(a-b)}{2}\right)^2\!\right)_+^{\frac{3}{2}}$.
\end{thm}
The Theorem \ref{thmweyl} can be extended to the pair operators $(H_0+V_1,\ H_0+V_2)$ (see Remark \ref{remarkweyl}).
\begin{rem}
The two formulae $(\ref{weyl})$ and  $(\ref{eqpremierterm})$ give in particular a Weyl type asymptotic of the counting function of the number of eigenvalues of $H_1$ between two values in the interval $]-mc^2,mc^2[\;$. In the case of scalar potential $v\,(v_+=v_-)$, this result was proved by Helffer-Robert \cite{BHDR} without the analytic assumption at infinity.
\end{rem}
To prove Theorem \ref{thmweyl} we construct, in Appendix A (see Theorem \ref{thmpropagator}), a parametrix at small time of the propagator of the Dirac equation in an external electro-magnetic field (see also Yajima \cite{KY} for a scalar electric potential cases).\\\\ 
As a direct result of the last theorem we deduce an upper bound $O(h^{-2})$ for the number of resonances close to non-critical energy levels in domains of width $h$  (see Proposition \ref{propmagh}) and a  Breit-Wigner approximation for the derivative of the spectral shift function $\xi(\la,h)$ (see Theorem \ref{thmbreit}).
\section{Distortion for the free Dirac operator}\label{secdistortion}
In this section, we start with the definition of
the deformation for the free Dirac operator by analytic distortion (in the spirit of Hunziker \cite{WH}) and we will calculate the essential spectrum for the distorded free Dirac operator. Here, $h$ does not play any role, and can be taken equal to $1$.
Let us now introduce the one-parameter family of unitary distortion 
$$U_{\theta}f(x)=J_{\phi_{\theta}(x)}^{\frac12}f(\phi_{\theta}(x)),\;
\,\, \theta\in \R,\,\,\, f\in (S(\R^3))^4,$$
where $\phi_{\theta}(x)=x+\theta g(x) $ and  $g:\R^3\longmapsto \R^3$ is a smooth function.
Let $J_{\phi_{\theta}(x)}\!\!=\!\det(I\!+\!\theta \nabla g(x))$ be the Jacobian of $\phi_{\theta}(x)$.\\ \\
We suppose that $g$ satisfies the assumption
$${\bf(A_g)}\left\{\begin{array}{l}
\mbox{(i)}\;\sup_{x\in \R^3}\|\nabla g(x)\|=M^{-1}<+\infty.\\\\
\mbox{(ii)}\;g(x)=0, \; in\; the\; compact\; B(0,R_0), \;(see\;(\ref{eqsector})) .\\\\
\mbox{(iii)}\;g(x)=x,\;\;  outside\; of\; a\; compact\; set\;K(\supset B(0,R_0)).
\end{array}\right.$$
\begin{lem}

For $\theta \in ]-M,M[$, $U_{\theta}$ can be extended as an unitary operator on ${\cal H}$.
\end{lem}
\Proof  Since $|\theta|<M$, we have $ \|\theta \nabla g(x)\|<1$, and
$$(\nabla\phi_{\theta}(x))^{-1}=(I+\theta
\nabla g(x))^{-1}=\sum_{n=0}^{\infty}(-1)^n(\theta)^n(\nabla g(x))^n.$$ 
The function $\phi_{\te}(x)$  is injective and  $\phi_{\te}(\R^3)=\R^3$, consequently
$\phi_{\theta}(x)$ is a diffeomorphism from  $\R^3$ to $\R^3$.
The inverse of $U_\te$ is given by $${U}_{\theta}^{-1}u=J_{\phi_{\theta}(x)}^{\frac{-1}{2}}u(\phi_{\theta}^{-1}(x)):\;(L^2(\R^3))^4 \longmapsto (L^2(\R^3))^4.$$
The lemma follows from the relations
\[U_{\theta}{U}_{\theta}^{-1}={U}_{\theta}^{-1}U_{\theta}=I_{(L^2(\R^3))^4}\; \mbox{and} \;\|
  U_{\theta}f\|_{{\cal H}}=\| {U}_{\theta}^{-1}f\|_{{\cal H}}=\| f\|,\,\,\, \forall f
  \in (L^2(\R^3))^4.\] \hfill{$\square$}
\begin{defi}
 We denote by ${\cal B}$, the space of functions  $f=(f_i)_{1\leq
  i\leq 4}$ such that  $f_i(x)$ has an analytic continuation in $C_{\epsilon,0}$ and  $\lim_{\!\!\!\begin{array}{c}
\scriptstyle |z|\rightarrow \infty \vspace{-2mm}\\
\scriptstyle z\in\, C_{\epsilon,0}
\end{array}}\!\!\!|z|^kf_i(z)=0$,\vspace{-2mm} for all $k\in \N$ and $\epsilon\in]0,1[$ (see (\ref{eqsector})).

\end{defi}
\begin{lem}The subspace
${\cal B}$ is dense in ${\cal H}$.
\end{lem}
\Proof The subspace ${\cal B}$ contains vectors of Hermite functions and the linear combinations of Hermite functions are dense in $L^2(\R^3)$.\hfill{$\square$}
\begin{prop}\label{propdense}
Let be $D_{\epsilon,M}=D_{\epsilon}\cap\{\theta \in \C; |\theta|<M\}$. We have the two assertions:
\begin{itemize}
\item[\emph{(i)}]For all $f\in {\cal B}$, $\theta\in D_{\epsilon,M}\longmapsto U_{\theta}f$ is analytic.
\item[\emph{(ii)}]For all $ \theta \in D_{\epsilon,M},U_{\theta}{\cal B}$ is dense in ${\cal H}$.
\end{itemize}
\end{prop}
\Proof In the order to prove (i), we show that 
$\theta \longmapsto\<U_{\theta}f,g\rangle$ is analytic for all $g \in {\cal H}$.
Let $R\gg1$ be such that $K\subset B(0,R)=\{x,|x|<R\}$.
 \begin{itemize}
 \item In $\{x,|x|<R\}$: since $f\in{\cal B}$, then $\theta\longmapsto\int_{|x|<R}J_{\phi_{\theta}}^{\frac 12}f(\phi_{\theta}(x))\overline{g}(x) dx$, is clearly analytic for all $g\in{\cal H}$.
\item In $\{x,|x|>R\}$: we have $ \;\; g(x)=x$, consequently  $\phi_{\theta}(x)=(1+\theta)x.$
We remark that
$$|\mbox{Im} (\phi_{\theta}(x))|=|\mbox{Im} (\theta x)|\leq
|\mbox{Im} (\theta)||x|=\frac{|\mbox{Im} (\theta)||\mbox{Re}(\phi_{\theta}(x))|}{|1+\mbox{Re}(\theta)|}.$$
If $|\theta|\leq r_{\epsilon}={\frac {\epsilon}{\sqrt{1+\epsilon^2}}},$  $0<\epsilon<1,$  then \[|\mbox{Im} (\phi_{\theta}(x))|\leq
\epsilon|\mbox{Re}(\phi_{\theta}(x))|.\] According to the definition of 
${\cal B}$ we have\[
|f(\phi_{\theta}(x))|\leq \frac{C_k}{|\phi_{\te}(x)|^k}\leq
\frac{C_{k,\epsilon}}{|x|^k},\;\; \forall |x|\geq R,\;\;  \forall\, \theta \in
D_{\epsilon,M},\;k\in\N\]
then, $\;\;\theta \longmapsto \int_{|x|\geq R}J_{\phi_{\theta}}^{\frac
  12}f(\phi_{\theta})\overline{g(x)} dx$ is analytic.
\end{itemize}
(ii) Let $h(x) \in (C_0^{\infty}(\R^3))^4$. We denote
$$h_k(x)=(\frac{k}{\pi})^{\frac{3}{2}}\int e^{-k(x-y-\theta
  g(y))^2}h(y)J_{\phi_{\theta}(y)}dy,$$
 which is clearly  in ${\cal B}.$\\
Using $(\frac{k}{\pi})^{\frac{3}{2}}\int e^{-k(x-y-\theta
  g(y))^2}J_{\phi_{\theta}(y)}dy=(\frac{k}{\pi})^{\frac{3}{2}}\int
e^{-k\,z^2}dz=1$, we get
$$h(x)-h_k(\phi_{\theta}(x))=(\frac{k}{\pi})^{\frac{3}{2}}\int
e^{-k(\phi_{\theta}(x)-\phi_{\theta}(y))^2}(h(x)-h(y))J_{\phi_{\theta}(y)}dy.$$
The last term tends to $0$ when $k\to+\infty$. Consequently, we have $$ h_k\circ\phi_{\theta}(x)\stackrel{k\to+\infty}{\To} h(x),\;\;\mbox{in}\;{\cal H}.$$\hfill{$\square$}
\begin{rem}
 One can always choose $g$ satisfying the assumption ${\bf(A_g)}$  with  $M>r_{\epsilon}={\frac {\epsilon}{\sqrt{1+\epsilon^2}}}$. In that case, we have  $D_{\epsilon,M}=D_{\epsilon}$. 
\end{rem}
\begin{lem}\label{lemmeH0}For $\te\in D_{\epsilon}$, we have 
$$H_{0,\theta}:=U_{\theta}H_0U_{\theta}^{-1}=\frac{1}{1+\theta}(-ic\sum_{j=1}^3\alpha_j\frac{\p}{\p x_j})+
\beta mc^2 + Q_{\theta}(x,\p_{x_j}),$$
where $Q_{\theta}(x,\p_{x_j})=\sum_{|\alpha|\leq
  1}a_{\alpha}(x,\theta)\p_{x_j}^{\alpha}$ is such that:
 \begin{itemize} 
\item [\emph{(i)}]$\theta\longmapsto a_{\alpha}(x,\theta)$ is an analytic function bounded by $O(\te)$.
\item [\emph{(ii)}]$x\longmapsto  a_{\alpha}(x,\theta) \in (C_0^{\infty}(\R^3))^4$.
\end{itemize}
In particular  $\theta\longmapsto H_{0,\theta}$ is an analytic family of type A of domain  $D(H_0)$ (see Kato \cite{TK}, for the definition of an analytic family of type A).
\end{lem}


\Proof
We denote $\p_j=\frac{\p}{\p x_j}$ and we calculate the term $U_{\theta}\p_j
U_{\theta}^{-1}$. \[U_{\theta}\p_j
U_{\theta}^{-1}f(x)=U_{\theta}\p_j\left(J_{\phi_{\theta}}^{\frac
  {-1}{2}}f(\phi_{\theta}^{-1}(x))\right)\;\;\;\;\;\;\;\;\;\;\;\;\;\;\;\;\;\;\;\;\;\;\;\;\;\;\;\;\;\;\;\;\;\;\;\;\;\;\;\;\;\;\;\;\;\;\]
\[\;\;\;\;\;\;\;\;\;\;\;\;\;\;\;\;= U_{\theta}\left(\p_j\left(J_{\phi_{\theta}}^{\frac
  {-1}{2}}\right).f\left(\phi_{\theta}^{-1}(x)\right)+J_{\phi_{\theta}}^{\frac
  {-1}{2}}\left(\p_jf(\phi_{\theta}^{-1}(x))\right)\right)\]
\[\;\;\;\;\;\;\;\;\;\;\;\;\;\;\;\;\;=-\frac {1}{2}
 J_{\phi_{\te}}^{-1}\p_jJ_{\phi_{\theta}(x)}f(x)+\sum_{k=1}^3\p_kf(x)\left(\p_j\phi_{{\te},k}^{-1}\right)(\phi_{\te}(x))
,\]
 with $\phi_{{\te}}^{-1}(x)=\left(\;\phi_{{\te},1}^{-1}(x),\;\phi_{{\te},2}^{-1}(x),\;\phi_{{\te},3}^{-1}(x)\right)$.\\
 We remark that   $\phi_{\te}^{-1}(x)=\frac{x}{1+\te}$   outside of the compact set $K$, then
\[\sum_{k=1}^3\p_kf(x)\left(\p_j\phi_{{\te},k}^{-1}\right)(\phi_{\te}(x))
=\frac{1}{1+\te}\p_jf(x),\;\; \mbox{outside of}\; K.\]
Let $\chi_K\in C_0^{\infty}(\R^3)$, $0\leq \chi_K\leq 1$, be equal to $1$ on $K$ and $0$ outside a compact set which contains $K$. We have
\[U_{\theta}\p_jU_{\theta}^{-1}f(x)=-\frac12 J_{\phi_{\te}}^{-1}\p_jJ_{\phi_{\theta}(x)}f(x)+\frac{1}{1+\te}\p_jf(x)(1-\chi_K)+\sum_{k=1}^3\p_kf(x)\left(\p_j\phi_{{\te},k}^{-1}\right)(\phi_{\te}(x))
\chi_K.\]
 \\Since $\p_jJ_{\phi_{\theta}(x)}$ has a compact support, then
 \begin{eqnarray}\label{eqQ}
 U_{\theta}\p_jU_{\theta}^{-1}=\frac{1}{1+\te}\p_j+q_{\theta}(x,\p_{x_j})
  \end{eqnarray} 
 with  $\;q_{\theta}(x,\p_{x_j})$ satisfying the hypothesis of  $Q_{\theta}(x,\p_{x_j})$.\\  
Now we just have to multiply ($\ref{eqQ}$) by  $-ic\al_j$, 
to sum on all values of $j$ and add  $\be mc^2$ to both hands. The estimate $a_{\alpha}(x,\theta)=O(\te)$ is clear using that $Q_0(x,\p_{x_j})=0$ and the analytic property of $\te\longmapsto a_{\alpha}(\cdot,\theta)$. \hfill{\hfill{$\square$}}

\begin{lem}\label{ess1} Let $P_{\theta}=\frac{1}{1+\theta}(-ic\sum_j\alpha_j\frac{\p}{\p x_j})+\beta
mc^2$. Then $$\sigma(P_{\theta})=\sigma_{ess}(P_{\theta})=\Gamma_\te=\{z\in \C; z=\pm
  c(\frac{\lambda}{(1+\theta)^2}+m^2c^2)^{\frac 12},\lambda \in [0,+\infty[\}.$$
\end{lem}
\Proof 
Let ${\cal F}$ be the Fourier transform  and 
\begin{eqnarray}K(\theta)={\cal F}\,P_{\theta}\,{\cal
  F}^{-1}&=&\frac{c}{1+\theta}\sum_j\alpha_j\xi_j+mc^2\beta\nonumber\\
  &=&\left(\begin{array}{cc}
mc^2I_2&(\frac{c}{1+\theta})(\sigma_1\xi_1+\sigma_2\xi_2+\sigma_3\xi_3)\\
(\frac{c}{1+\theta})(\sigma_1\xi_1+\sigma_2\xi_2+\sigma_3\xi_3)&-mc^2I_2 
 \end{array}\right)\nonumber
 \end{eqnarray}  
   where
$\xi=(\xi_1,\xi_2,\xi_3)\in \R^3$ and  $\alpha_j\xi_j$ is the multiplication operator  by  the $4\!\times\!4$ matrix $\alpha_j\xi_j$. Here $\sigma_1,\sigma_2,\sigma_3$ are the $2\times2$ Pauli matrices. The spectrum of $P_{\te}$ coincide with the spectrum of the multiplication operator $K(\te)$. We easily prove that 
$$\si(K(\te))=\si_{ess}(K(\te))=\Gamma_\te=\{z\in \C; z=\pm
  c(\frac{\lambda}{(1+\theta)^2}+m^2c^2)^{\frac 12},\lambda \in [0,+\infty[\},$$
and we deduce the lemma.\hfill{$\square$}\\
 The principal branch of the square root function is holomorphic on the set $\C\setminus]-\infty,0]$. Let $S_{D_{\epsilon}}=\{z=\frac{\lambda}{(1+\theta)^2}+m^2c^2,\;\te\in\,D_{\epsilon},\;\lambda \in [0,+\infty[\}$. Since, $$S_{D_{\epsilon}}\subset ]0,+\infty[e^{]-\frac{\pi}{2},\frac{\pi}{2}[},$$ the square root  $z\longmapsto z^{\frac12}$ is holomorphic on $S_{D_{\epsilon}}$.
\begin{lem}\label{ess2}For $H_{0,\theta},\,P_{\theta}$, defined as above, we have 
$\sigma_{ess}(H_{0,\theta})=\sigma_{ess}(P_{\theta})$.

\end{lem}
\Proof We want to use Kato's theorem \cite[Th.4.5.35]{TK}.
For $\lambda\gg1,\;\lambda \in \R$ and $Q_{\te}$ defined in Lemma \ref{lemmeH0}, we have
$$(H_{0,\theta}-i\lambda)=(1+Q_{\theta}(P_{\theta}-i\lambda)^{-1})(P_{\theta}-i\lambda).$$ Since
$(P_{\theta}-i\lambda)^{-1}\in {\cal L}({\cal H},(H^1)^4)$ and
$Q_{\theta}(P_{\theta}-i\lambda)^{-1}=O(\frac{\te}{\la})$, we obtain that $i\lambda \in \rho
(H_{0,\theta})=\C\setminus \sigma(H_{0,\theta}).$
To apply the Kato Theorem, it is enough to show that
\begin{eqnarray}\label{eqcompact}
(H_{0,\theta}-i\lambda)^{-1}-(P_{\theta}-i\lambda)^{-1} \mbox{ is compact}.
\end{eqnarray}
Using the resolvent equation, we have
$$(H_{0,\theta}-i\lambda)^{-1}-(P_{\theta}-i\lambda)^{-1}=(H_{0,\theta}-i\lambda)^{-1}Q_{\theta}(P_{\theta}-i\lambda)^{-1}.$$
with $Q_{\theta}(x,\p_{x_j})=\sum_{|\alpha|\leq
  1}a_{\alpha}(x,\theta)\p_{x_j}^{\alpha}$ compactly supported. Since the operaror $(H_{0,\theta}-i\lambda)^{-1}Q_{\theta}$ is bounded 
  and ${\1}_{ \mbox{supp}(Q_{\te})}(P_{\theta}-i\lambda)^{-1} $ is compact,
 the assertion (\ref{eqcompact}) holds.
\hfill{$\square$}
\section{Definition of resonances }\label{secdefinition}
In this section we distort the perturbed Dirac operator $H=H_0+V$, where the potential $V$ satisfies the assumption ${\bf(A_V)}$ and we define the resonances for the semi-classical Dirac operator.\\
 The distorted Dirac operator is denoted by $$H_{\theta}=U_{\theta}H_0U_{\theta}^{-1}+U_{\theta}VU_{\theta}^{-1}=H_{0,\theta}+V(\phi_{\te}(x)).$$ 

\begin{prop}\label{propana}
We suppose that the potential $V$ satisfies the assumption ${\bf(A_V)}$
, then
\begin{itemize}
\item[\emph{(i)}]$\theta\in D_{\epsilon}\longmapsto H_{\theta}=H_{0,\theta}+
V(\phi_{\theta}(x))$ is an analytic family of type A.
\item[\emph{(ii)}]$\sigma_{ess}(H_{\theta})=\Gamma_\te$.
\end{itemize}
\end{prop}
\Proof
The assertion (i) is clear since $H_{0,\theta}$ is an analytic family of type A and $V$ satisfies the assumption ${\bf (A_V)}$.\\
Now, we prove (ii) as in the proof of Lemma \ref{ess2}. 
For $\la\gg1$, $i\la\in \rho(H_{\te})$ and 
\begin{eqnarray}\label{eqcpt}(H_{\te}-i\la)^{-1}-(H_{0,\te}-i\la)^{-1}= (H_{\te}-i\la)^{-1}V(\phi_{\te}(x))(H_{0,\te}-i\la)^{-1}.
\end{eqnarray}
Since the operator $V(\phi_{\te}(x))(H_{0,\te}-i\la)^{-1}$ is compact (see (${\bf A_V}$)) and the resolvent  $(H_{\te}-i\la)^{-1}$  is bounded, the difference $(H_{\te}-i\la)^{-1}-(H_{0,\te}-i\la)^{-1}$ is compact.
According to Kato's theorem \cite[Theorem.4.5.35]{TK} and to the lemmas $\ref{ess1},\ref{ess2}$, we obtain (ii).\hfill{$\square$}

We denote $$\Sigma=\{z\in\C;\;\mbox{Im} (z)\geq0,\;\mbox{Re}(z)>-mc^2\}\cup\{z\in\C;\;\mbox{Im}\, (z)\leq0,\;\mbox{Re}(z)<mc^2\}\setminus \si(H).$$
\begin{thm}\label{thmdef}With the notations used above and taking $\theta_0 \in D_{\epsilon}^{+}=D_{\epsilon}\cap
\{\emph{Im} (\theta)\geq0\}$, we have:
\begin{itemize}
\item[\emph{(i)}]For all $ f,g \in {\cal B}$, the function: $z\in\Sigma\longmapsto
M_{f,g}(z)=\<(z-H)^{-1}f,g\rangle$  has a meromorphic extension on 
$S_{\te_0}$.
\item[\emph{(ii)}]The poles of $M_{f,g}(z)$  are the eigenvalues of 
 $H_{\theta_0}$. 
\item[\emph{(iii)}]These poles are independent of the family  $U_{\theta_0}$.
\item[\emph{(iv)}]$\sigma_d(H_{\theta_0})\cap\Sigma=\emptyset,$
where $\si_d(H_{\te_0})$ is the discrete spectrum of the operator $H_{\te_0}$.
\end{itemize}
\end{thm}
\Proof
(i) Since $U_{\theta}$ is unitary for  $\theta \in \R$,
\begin{eqnarray*}M_{f,g}(z)=\<(z-H)^{-1}f,g\rangle
=\<(z-H_{\theta})^{-1}U_{\theta}f,U_{\theta}g\rangle.\end{eqnarray*}
 We denote
\begin{eqnarray}\label{eqdefM}M_{f,g,\theta}(z)=\<(z-H_{\theta})^{-1}U_{{\theta}}f,U_{\bar{\theta}}g\rangle,\;\;
\mbox{for}\; \theta \in D_{\epsilon}.\end{eqnarray}
According to (i) of Proposition \ref{propana} and to the definition of $U_{\te}$, the functions 
$\theta\mapsto(z-H_{\theta})^{-1}$,
  $\theta\mapsto U_{\theta}f$ and $\theta \mapsto
\<\psi,U_{\bar{\theta}}g\rangle$ are analytic on $ D_{\epsilon}$ for all $ \psi \in {\cal H}$ and any $z\in\Sigma$.

Thus, for $z\in\Sigma$, the function $\theta \mapsto M_{f,g,\theta}(z)$
is analytic on $D_{\epsilon}$.
Since $M_{f,g,\theta}(z)$ is independent of  $\te$ on the real axis and according to the uniqueness of the prolongation, its extension is independent of $\te $.

Now, we fix  $\theta_0 \in D_{\epsilon}^{+}$. 
 Since $S_{\te_0} \cap \si_{ess}(H_{\theta_0})=\vide$, the function $z\in\Sigma\longmapsto M_{f,g,\theta_0}(z)$ has a meromorphic extension in $S_{\te_0}$.\\ 

(ii) First, let $z\in S_{\te_0}$  be a pole of 
$M_{f,g}(z)$ which is equal to  $M_{f,g,\theta_0}(z)$ for $\theta_0\in D_{\epsilon}^{+}$. Then $z\in \sigma_d(H_{\theta_0})\cap S_{\te_0}$ (see proof of (i)).

Now, let $e\in \sigma_d(H_{\theta_0})\cap S_{\te_0}$.
There exists  $u\in {\cal H}$ such that $\| u \Vert =1$ and $H_{\theta_0}
u=eu$.
Let $\gamma$ be a small disk centered at $e$ such that  ${\gamma}^{^{\!\!\!\!\circ}}\cap
\sigma (H_{\theta_0})=\{e\}$ and $\Gamma$ be the positively oriented boundary of $\gamma$.\\ 
Let us introduce the projector
$$\Pi=\frac{1}{2i\pi}\int_{\Gamma}(z-H_{\te_0})^{-1}dz;\;\; \Pi u=u.$$ Since
${\overline{U_{\theta_0}{\cal B}}={\cal H}}={\overline{U_{\overline{\theta}_0}{\;\cal B}}}$ (see Proposition \ref{propdense}),
 there exist $f_n,\;g_n\in {\cal B}$ such that \[\vert
u-U_{\theta_0}f_n\vert\leq \frac1n\; \;\mbox{and}\; \;\vert u-U_{{\bar\theta}_0}g_n \vert \leq \frac1n,\;\;\;n\in\N.\] 
Therefore, as $n$ goes to infinity, we have
\begin{eqnarray*}\frac{1}{2i\pi}\!\int_{\Gamma}\!\!
\<(z-H_{\theta_0})^{-1}U_{\theta_0}f_n,U_{\bar{\theta_0}}g_n\rangle dz&=&\frac{1}{2i\pi}\!\int_{\Gamma}\!\!\<(z-H_{\theta_0})^{-1}u,u\rangle dz
+o(1)
\\&=&\<\Pi u,u\rangle+o(1)\\&=&\| u\|^2+o(1)\\&=&1+o(1),\end{eqnarray*} 
and then,
\[\frac{1}{2i\pi}\int_{\Gamma}\<(z-H)^{-1}f_n,g_n\rangle dz=1+o(1).\]
So that $M_{f_n,g_n}(z)$ admits $e$ as a pole  in  $\gamma$.

The assertion (iii) results from  (ii) because $M_{f,g}(z)$ is independent of $U_{\theta}$.

(iv) If there exists  $z\in \sigma_d(H_{\theta_0})\cap\Sigma$ 
, then $ z$ is a pole of 
$\<(z-H)^{-1}f,g\rangle,\;\;$for $f,g\in {\cal B}$, but $\<(z-H)^{-1}f,g\rangle$ is analytic on this domain. We conclude that such $z$ does not exist.\hfill{$\square$}  
\begin{rem}\label{remark6}
\begin{itemize}
  \item[(i)] It results from $(ii)$ of Theorem \ref{thmdef} that for all $\te\in D_{\epsilon}^{+}$,
the discrete spectrum $\sigma_d(H)$ is a subset of $\sigma_d(H_{\theta})$.\\
\item[(ii)] The previous theorem justifies the definition of the resonances (Definition \ref{defres}) and using Lemma \ref{ess1}, $H_0$ has no resonances.
\end{itemize}
\end{rem}
\begin{rem}
  If $\te\in D_{\epsilon}$, then its conjugate $\bar \te \in D_{\epsilon}$. Repeating arguments of the proof of Theorem \ref{thmdef}, we have
\begin{itemize}
  \item[(i)] The function $\theta \longmapsto M_{f,g,\bar\theta}(z)$ has a analytic extension  for $\te\in D_{\epsilon}$.
  \item[(ii)] The function $z\in\bar\Sigma\longmapsto M_{f,g,\bar\theta}(z)$  has a meromorphic extension on $S_{\bar\te_0}$,
  \end{itemize}
  where $$\bar\Sigma=\{z\in\C;\emph{Im} (z)\geq0,\;\emph{Re}(z)<mc^2\}\cup\{z\in\C;\emph{Im} (z)\leq0,\;\emph{Re}(z)>-mc^2\}\setminus \si(H),$$and $S_{\bar\te_0}$ is the symmetric of $S_{\te_0} $ with respect to the real axis
  \[S_{\bar\te_0}=\{ z\in\bigcup_{
\substack{\te\in D_{\epsilon}^+}}\!\!\Gamma_{\bar\te}\,;\;\;arg(1+\te)<\mbox{arg}(1+\te_0),\;\;\frac{1}{|1+\te|}<\frac{1}{|1+\te_0|}\;\}.\vspace{-3mm}\]\vspace{-2mm}
Consequently, we obtain (see Theorem \ref{thmdef}):\\\\
$1)$ The poles of $M_{f,g}(z)$ in $S_{\bar\te_0}$ are the eigenvalues of 
 $H_{\bar\theta_0}$. \\
$2)$ These poles are independent of the family  $U_{\bar\theta_0}$.\\
$3)$ $\sigma_d(H_{\bar\theta_0})\cap 
\bar\Sigma=\emptyset$.\\
The assertions $3)$ and \emph{(ii)} prove that the operator $H_{\bar \theta}$ has only discrete spectrum in $S_{\bar\te_0}$.  
\end{rem}
{\bf Proof of the Proposition \ref{propimmaginaire}.}
$\label{proof1}$ We consider the anti-linear application on $\HH$, $$C:\psi\longmapsto U_c{\overline\psi}=i\be\al_2{\overline\psi}.$$
 Then, we have
 \begin{eqnarray*}C\;H_{\te}\;C^{-1}&=&-H_{0,\bar \te}+C V\circ\phi_{\te}(x) C^{-1}\\
 &=&-\left(H_{0,\bar \te}- U_c{\overline {V\circ\phi_{\te}(x)}}U_c^{-1}\right).\end{eqnarray*}
$\mbox{Using that } {V}\mbox{ is analytic,  } \mbox{ we get}\; \overline {V\circ\phi_{\te}(x)}={\overline {V}\circ\phi_{\bar \te}(x)}$.
Then,
\begin{eqnarray*}C\;H_{\te}\;C^{-1}&=&-\left(H_{0,\bar \te}- U_c{\overline {V}\circ\phi_{\bar \te}(x)}U_c^{-1}\right)\\&=&-\left(H_{0,\bar \te}- U_cU_{\bar \te}{\overline {V}(x)}U_{\bar \te}^{-1}U_c^{-1}\right).\end{eqnarray*}
We recall that $U_{\theta}f(x)=J_{\phi_{\theta}(x)}^{\frac12}f(\phi_{\theta}(x))$. Since $U_cU_{\te}=U_{\te}U_c$, we obtain
\begin{eqnarray*}C\;H_{\te}\;C^{-1}&=&-\left(H_{0,\bar \te}- U_{\bar \te}U_c{\overline {V}(x)}U_c^{-1}U_{\bar \te}^{-1}\right)\\&=&-U_{\bar\te}\left(H_{0}-U_c{\overline {V}(x)}U_c^{-1}\right)U_{\bar\te}^{-1}=-H_{\bar\te}^-.\end{eqnarray*}
Consequently,$\;\;\;C\;(H_{\te}-z)\;C^{-1}=-(H_{\bar\te}^-+\bar z),$ and the property follows.\hfill{$\square$}\\\\
{\bf Proof of the Proposition \ref{propreelle}.}$\label{proof2}$ \\
By definition of $H_{\te}$, we have
\begin{eqnarray}\label{eqprop2}
\overline{H_{\te}-z}=U_{\bar\te}\left({\overline H}_0+{\overline V}-\bar z\right)U_{\bar\te}^{-1}.
 \end{eqnarray}
Using that  $\bar\al_1=\al_1,\;\bar\al_2=-\al_2,\;\bar\al_3=\al_3,\;\bar\be=\be$, we find 
\[{\overline H}_0=ic\sum_{j=1}^{3}\overline{\alpha}_j \frac{\p }{\p x_j} +
\beta mc^2=-ic\sum_{j=1}^{3}\alpha_j' \frac{\p }{\p x_j} +
\beta' mc^2\;\;\]and\[\;\; {\overline H}_0+{\overline V}={\overline H}.\]
 Using  the last relation and the equation (\ref{eqprop2}), we obtain Proposition \ref{propreelle}.\hfill{$\square$}\\\\
 Finally, the study of resonances in a domain of the complex plan $\C$ is reduced to study the resonances in  $\Om\cap\{z\in \C,\ \mbox{Im} (z)<0\}$, with $\Om$ satisfying the assumption ${\bf(A_{\Om}^+)}$ (see Fig.1).

\section{ Upper bound for the number of resonances}\label{majoration}
In this section, we establish an upper bound of the number of resonances in a compact domain $\Om$. For this purpose we construct an operator ${\widehat H}_{\te}:D(H_0)\to \HH$ with some properties (see Proposition \ref{propconst}).
According to the Section \ref{secdefinition}, it is sufficient to treat the case where $\Om$ satisfies the assumption ${\bf(A_{\Om}^+)}.$

We shall use the theory of $h$-pseudo-differential operator (see \cite{MDJS}, \cite{DR1}). 
Let $m$ be an order function on $\R^{2n}$ (i.e. there are $C_0,N_0>0$, such that $m(x)\!\leq\! C_0 \langle x-y\rangle^{N_0}m(y)$).
The space ${\cal S}^p(m)$ is the set of $a(x,\xi;h)\in C^\infty(\R^{2n})\otimes\C^4$ such that for every $\al\in\N^{2n}$, there exists $C_\al>0$, such that $$\|\nabla_{x,\xi}^\al a(x,\xi;h)\|\leq C_\al m(x,\xi)h^{-p}.$$
For a symbol $a(x,\xi;h)$, we define the Weyl quantization, $a^w(x,h\nabla_x;h):=\Op(a)$ by
$$\Op(a)u(x)=\frac{1}{(2\pi h)^n}\int_{\R^n}\!\!\int_{\R^{n}}e^{ih^{-1}(x-y)\cdot\xi}a(\frac{x+y}{2},\xi; h)u(y)dyd\xi,$$
where $u(x)$ is in the Schwarz 
 space.
\subsection{Construction of ${\widehat H}_{\te}$}\label{Appendix A}We follow the approach of Sj\"ostrand \cite{JS2}.
Let $\Om$ be a complex domain satisfying the assumption ${\bf(A_{\Om}^+)}$ and  $\psi\in C_0^{\infty}(\R^3)$ be such that $\psi(x)\geq0,\;\psi(x)=1$ if $|x|\leq1$ and $\psi(x)=0$ if $|x|\geq2$. 
We recall the notations of Section $\ref{secdistortion}$: $\phi_{\te}(x)=x+\te g(x)$ with $g(x)=0$ in the compact set $B(0,R_0)\subset K$ and $g(x)=x$ outside $K\subset B(0,\al_0)$ where $\al_0>0$ is sufficiently large.

Using  Lemma $\ref{lemmeH0}$, 
the semi-classical principal symbol of $H_{\te}$ is given by:
$$h_{\te}(x,\xi)=\al\cdot\zeta_{\te}(x,\xi)+mc^2\be+V(\phi_{\te}(x))-c\frac12\sum_{j=1}^3\al_j\cdot J_{\phi_{\te}(x)}^{-1}\p_jJ_{\phi_{\te}(x)},$$ 
 with 
$$\zeta_{\te}(x,\xi)=\left(\zeta_{\te,1}(x,\xi),\;\zeta_{\te,2}(x,\xi),\;\zeta_{\te,3}(x,\xi)\right)\;\;\mbox{and}\;\;\zeta_{\te,j}(x,\xi)=c\sum_{k=1}^3\xi_k\left(\p_j\phi_{\te,k}^{-1}\right)(\phi_{\te}(x)).$$
For all $(x,\xi)$, the matrix $M=\al\cdot\zeta_{\te}(x,\xi)+mc^2\be$, has two eigenvalues 
$$\la_{\te}^{\pm}=\pm\sqrt{\zeta_{\te}(x,\xi)^2+m^2c^4}.$$ 
Consequently, there exists an invertible matrix $U$  such that $$U^{-1}MU=d_{\te}:=\left(\begin{array}{cc}
\la_{\te}^+I_2&0\\
0&\la_{\te}^-I_2
\end{array}\right),$$where
\[U=\left(\begin{array}{cc}
I_2&\frac{-1}{\la_{\te}^++mc^2}\si\cdot\zeta_{\te}(x,\xi)\\
\frac{1}{\la_{\te}^++mc^2}\si\cdot\zeta_{\te}(x,\xi)&I_2\\
\end{array}\right),\]
with $\si=(\si_1,\si_2,\si_3)$  and $(\si_j)_{1\leq j\leq3}$ the $2\!\times\!2$ Pauli matrices.\\
One can easily prove that the norms of $U,\,U^{-1}$ and their derivatives are bounded in the following way:
\begin{eqnarray}\label{eqmajU}
\|\p_x^\al \p_\xi^\be U(x,\xi)\|&<&C\langle\xi\rangle^{-\be},\nonumber\\
\|\p_x^\al \p_\xi^\be U^{-1}(x,\xi)\|&<&C\langle\xi\rangle^{-\be},\;\;\;\forall\al,\,\be\in\N.
\end{eqnarray}
 Applying $U^{-1}$(resp. $U$) on the left (resp. on the right) of $h_{\te}(x,\xi)$, we obtain
 \[U^{-1}h_{\te}(x,\xi)U=d_{\te}+{\Ti V}_{\te}(x,\xi),\]
 where ${\Ti V}_{\te}(x,\xi)=U^{-1}\left(V(\phi_{\te}(x))-c\frac12\sum_{j=1}^3\al_j\cdot J_{\phi_{\te}(x)}^{-1}\p_jJ_{\phi_{\te}(x)}\right)U$. Since the term $\p_j\!\left(J_{\phi_{\te}(x)}\right)\!$ is compactly supported $j\!=\!1,2,3,$ (see Lemma $\ref{lemmeH0}$) and the matrix $V(\phi_{\te}(x)),\,U,\;U^{-1},$ and their derivatives are uniformly bounded, then ${\Ti V}_{\te}(x,\xi)$ and their derivatives are uniformly bounded.\\\\
In order to construct  ${\widehat H}_{\te}$, we introduce an intermediate 
 function $f(x,\xi)$:\\\\
We denote $|\Om|$ the diameter of $\Om$. Let us choose $\be_0>0$ and $C_0>0,$ sufficiently large such that   \begin{eqnarray}\label{eq**}\forall\,\xi\, \in\, \R^3;\;\; \sup_{x\in\,\R^3}\|{\Ti V}_{\te}(x,\xi)\|+|\Om|&\!\!\!\leq&\!\!\! \frac12\,|\la_{\te}^{\pm}-iC_0\psi(\frac{\xi}{\be_0})|\\ \nonumber
&\!\!\!=&\!\!\!\frac12\sqrt{\!\left(\!\mbox{Re}(\la_{\te}^{\pm})\right)^2\!+\!\left(\mbox{Im} (\la_{\te}^{\pm})\!-C_0\psi(\frac{\xi}{\be_0})\!\!\right)^2}.
\end{eqnarray} 
We prove $(\ref{eq**})$ considering the two cases:
\begin{itemize}
\item For $|\xi|>\be_0$, with $\be_0>0,$ sufficiently large we have
$$ \sup_{x\in\,\R^3}\|{\Ti V}_{\te}(x,\xi)\|+|\Om|\leq \frac12\,|\mbox{Re}(\la_{\te}^{\pm})|.$$ 
\item For $ |\xi|\leq \be_0$, since $\la_\te^\pm$ is bounded, we choose $C_0>0,$ sufficiently large such that
$$ \sup_{x\in\,\R^3}\|{\Ti V}_{\te}(x,\xi)\|+|\Om|\leq \frac12\,|\mbox{Im} (\la_{\te}^{\pm})-C_0|.$$ 
\end{itemize}
For $|x|>\al_0>0,$ sufficiently large we have $\zeta_{\te,j}(x,\xi)=\frac{c\xi_j}{1+\te}$ and $\;\la_{\te}^{\pm}=\pm c\sqrt{\frac{\xi^2}{(1+\te)^2}+m^2c^2}$. Since the domain $\Om$ satisfies the assumption ${\bf (A_{\Om}^+)}$, we have $$\mbox{min}\{\mbox{dist}(\overline\Om,\la_{\te}^+),\;\mbox{dist}(\overline\Om,\la_{\te}^-)\}\neq0,$$ hence we can choose $\al_0>0$ sufficiently large such that 
\begin{eqnarray}\label{eqVtilde}
\forall |x|>\al_0,\;\;\;\;\|{\Ti V}_{\te}(x,\xi)\|\leq\frac12 \;\mbox{dist}(\overline\Om,\la_{\te}^{\pm}):=\frac12 \;\mbox{min}\{\mbox{dist}(\overline\Om,\la_{\te}^+),\;\mbox{dist}(\overline\Om,\la_{\te}^-)\}.
\end{eqnarray}
Now, we define $f(x,\xi)$ in the following way:
\begin{eqnarray}\label{eqdeff}
f(x,\xi)=C_0\psi(\frac{x}{\al_0})\psi(\frac{\xi}{\be_0}).
\end{eqnarray}
 
\begin{lem}\label{Htilde}
$\;$The matrix $h_{\te}(x,\xi)-if(x,\xi)-z$ is invertible for all $z\in\Om$ and  satisfies 
\begin{eqnarray}\label{eqnormsymbol}
\|\p_x^\al \p_{\xi}^{\be}\left(h_{\te}(x,\xi)-if(x,\xi)-z\right)^{-1}\|< C\<\xi\rangle^{-1-\be},\;\;\;\forall\al,\,\be\in\N.
\end{eqnarray}
\end{lem}
\Proof
Applying $U^{-1}$(resp. $U$) on the left (resp. on the right) of $h_{\te}(x,\xi)-if(x,\xi)-z$, we obtain
$$U^{-1}\left(h_{\te}(x,\xi)-if(x,\xi)-z\right)U=d_{\te}-if(x,\xi)-z+{\Ti V}_{\te}(x,\xi).$$
1) Let us prove that the symbol 
 ${{\cal \si}}:= d_{\te}-if(x,\xi)-z+{\Ti V}_{\te}(x,\xi)$ is invertible.
\begin{itemize}
\item For $|x|\leq\al_0$,
\begin{eqnarray}
{ {\cal \si}}&=&(d_{\te}-iC_0\psi(\frac{\xi}{\be_0}))\left(I_4+\left(d_{\te}-iC_0\psi(\frac{\xi}{\be_0})\right)^{-1}\left({\Ti V}_{\te}(x,\xi)-z\right)\right).\nonumber
\end{eqnarray}
 According to $(\ref{eq**})$, we have $$\|\left(d_{\te}-iC_0\psi(\frac{\xi}{\be_0})\right)^{-1}\left({\Ti V}_{\te}(x,\xi)-z\right)\|<\frac12,$$thus ${ {\cal \si}}$ is invertible and satisfies
 \begin{eqnarray}\label{spm}
 \| { {\cal \si}}^{-1}\|<2\Big \|\left(d_{\te}-iC_0\psi(\frac{\xi}{\be_0})\right)^{-1}\Big\|<C\<\xi\rangle^{-1}.
 \end{eqnarray}
\item For $|x|>\al_0,$ we have $\;\la_{\te}^{\pm}=\pm c\sqrt{\frac{\xi^2}{(1+\te)^2}+m^2c^2}$.
 Since $f(x,\xi)\geq0$, we have 
$$|\la_{\te}^+-(z+if(x,\xi))|>\mbox{dist}(\overline\Om,\la_{\te}^{\pm})>C\<\xi\rangle>0,$$ 
and
$$|\mbox{Re}\left(\la_{\te}^--(z+if(x,\xi))\right)|=|\mbox{Re}\left(\la_{\te}^--z\right)|>\mbox{dist}(\overline\Om,\la_{\te}^{\pm})>C\<\xi\rangle>0.$$
Since 
\begin{eqnarray}
{ {\cal \si}}&=&(d_{\te}-if(x,\xi)-z)\left(I_4+\left(d_{\te}-if(x,\xi)-z\right)^{-1}{\Ti V}_{\te}(x,\xi)\right),\nonumber
\end{eqnarray}
and $$\|\left(d_{\te}-if(x,\xi)-z\right)^{-1}{\Ti V}_{\te}(x,\xi)\|<\frac12,\;\;\mbox{(see (\ref{eqVtilde}))},$$
the matrix ${ {\cal \si}}$ is invertible and 
\begin{eqnarray}\label{spm1}
 \|{ {\cal \si}}^{-1}\|<2\|\left(d_{\te}-if(x,\xi)-z\right)^{-1}\|<C\<\xi\rangle^{-1}.
\end{eqnarray}
\end{itemize}
2) According to 1), the matrix 
$U^{-1}(h_{\te}(x,\xi)-if(x,\xi)-z)U$ 
is invertible. From (\ref{spm}), (\ref{spm1}) and (\ref{eqmajU}), we deduce 
that the matrix $h_{\te}(x,\xi)-if(x,\xi)-z$ is invertible and
 \begin{eqnarray}\label{eqbe=0}\|\left(h_{\te}(x,\xi)-if(x,\xi)-z\right)^{-1}\|&=&\|U^{-1}\left(d_{\te}-if(x,\xi)I_4-z+{\Ti V}_{\te}(x,\xi)\right)^{-1}U\| \nonumber\\ \nonumber
 &\leq& \|U\|\|U^{-1}\|\|\left(d_{\te}-if(x,\xi)I_4-z+{\Ti V}_{\te}(x,\xi)\right)^{-1}\|\\
 &<&C\<\xi\rangle^{-1}.
 \end{eqnarray}
This gives (\ref{eqnormsymbol}) for $\al=\be=0$. Using (\ref{eqmajU}) and (\ref{eqbe=0}) we obtain (\ref{eqnormsymbol}) for $(\al,\be)\in\N^2$ by induction.

 \hfill{$\square$}\\
 We denote ${\Ti H}_{\te}=H_{\te}+{\Ti T},$ with ${\Ti T}=\Op(-if(x,\xi))$, where $f(x,\xi)$ is defined in (\ref{eqdeff}).
  It's clear that the semi-classical principal symbol of  $({\Ti H}_{\te}-z)$ is  $${\cal \si}_{\Ti H_{\te}}:=h_{\te}(x,\xi)-if(x,\xi)-z.$$
\begin{prop}
If $h>0$ is small enough, the operator $(z-{\Ti H}_{\te})$ 
is invertible for every $z\in \Om$ and, for every $N\in \N$ its inverse satisfies:
$$\;\;\;\;\;\;\;\;\;\hspace{5mm}(z-{\Ti H}_{\te})^{-1}=O_N(1): D(H^N)\longmapsto D(H^{N+1}),\;$$ uniformly for $z\in \Om $. Here $D(H^N)$ designates the domain of $H^N$ $\mbox{with the convention}\; \ D(H^0)=\HH$.
\end{prop}
\Proof
Let us prove that the operator $(z-{\Ti H}_{\te})$ is a Fredholm operator of index $0$. We have
\begin{eqnarray}
(z-{\Ti H}_{\te})(z-H_{0,\te})^{-1}&=&(z-H_{0,\te}+H_{0,\te}-{\Ti H}_{\te})(z-H_{0,\te})^{-1}\nonumber\\\nonumber
    &=&I-({\Ti T}+V(\phi_{\te}(x)))(z-H_{0,\te})^{-1}.
\end{eqnarray}
Since the right-hand side is a perturbation of the identity by a compact operator and$$(z\!-\!H_{0,\te})^{-1}:(L^2(\R^3))^4\mapsto D(H)\;\mbox{ is invertible},$$ the
 operator $(z-{\Ti H}_{\te})$ is Fredholm of index $0$.
 Consequently, it is enough to show that 
 \begin{eqnarray}{\label{eqDNDN}}
 \|u\|_{D(H^{N+1})}^2\leq C \|(z-{\Ti H}_{\te})u\|_{D(H^N)}^2,\;\;\;\mbox{for }\; u\in D(H^{N+1}).
 \end{eqnarray}
 According to Lemma $\ref{Htilde}$, the symbol $q_0={\cal \si}_{\Ti H_{\te}}^{-1}$ is well defined and satisfies
\[\|\p_x^\al \p_{\xi}^{\be}q_0\|< C\<\xi\rangle^{-1-\be}.\]
Moreover, having
\[\|\p_x^\al \p_{\xi}^{\be}({\cal \si}_{\Ti H_{\te}})\|< C\<\xi\rangle^{+1-\be},\]
 the composition theorem of  $h$-pseudo-differential operators implies
$$\Op(q_0)\Op({\cal \si}_{\Ti H_{\te}})=\Op(r)$$
where $(r-1)$ is in the space of symbols ${\cal S}^{0}(h)$.
 In particular the operator $$\Op(r):\;D(H^{N+1})\longmapsto D(H^{N+1}),\;\;\;\forall N\in\N,$$ is invertible for $h$ small enough, then (\ref{eqDNDN}) follows. Therefore the operator $({\Ti H}_{\te}-z)$ is also invertible and we have $$(z-{\Ti H}_{\te})^{-1}=O_N(1): D(H^N)\longmapsto D(H^{N+1}).$$
 \hfill{$\square$}

\begin{prop}\label{propconst}
There exists ${\widehat H}_{\te}:D(H)\longmapsto \HH$, with the following properties.\\ 
The difference $K:={\widehat H}_{\te} -H_{\te}$ is of finite rank $O(h^{-3})$, has a compact support in the sense that $K=\chi_1 K\chi_1$ for some $\chi_1 \in C_0^{\infty}(\R^3)$ and 
  \[K=O(1):D(H^N)\longmapsto D(H^M)\;\forall N,M \in \N.\] Moreover, for every  $N\in \N$, we have$$({\widehat H}_{\te}-z)^{-1}=O(1):D(H^N)\longmapsto D(H^{N+1}),\; $$uniformly for $z\in {\overline \Om}$. 
\end{prop}
\Proof (We again use all the previous notations) 
We define $${\widehat H}_{\te}:=H_{\te} +{\chi}_1T\chi_1={\Ti H}_{\te}+{\chi}_1T{\chi}_1-{\Ti T},$$
with $\chi_1(x)=\psi(\frac{x}{2\al_0})$ 
and
$$T:=\chi(-h^2\De+x^2){\Ti T}=\chi(-h^2\De+x^2)\Op(-if(x,\xi))$$
where $\;\chi \in C_0^{\infty}(\R)$ is such that: 
$$\chi(\xi^2+x^2)=1\mbox{ on the support of }f(x,\xi) \ (\mbox{see}\ (\ref{eqdeff})).$$
By functional calculus (see \cite{MDJS}), we can prove that 
\begin{eqnarray}\label{eqtauinfinity}
{\Ti H}_{\te}-{\widehat H}_{\te}={\Ti T}-{\chi}_1T\chi_1=O(h^{\infty}):D(H^N)\longmapsto D(H^M),\;\forall \;M,\ N\in \N.
\end{eqnarray}
The last lemma, the formula (\ref{eqtauinfinity}) and  
\[({\widehat H}_{\te}-z)^{-1}=({\Ti H}_{\te}-z)^{-1}\left(I+({\widehat H}_{\te}-{\Ti H}_{\te})({\Ti H}_{\te}-z)^{-1}\right)^{-1}\]
yield for all N $\in \N$
\[({\widehat H}_{\te}-z)^{-1}=O(1): D(H^N)\longmapsto D(H^{N+1}).\;\;\;\]
According to the fact that $\chi(-h^2\De+x^2)$ is of finite rank  $O(h^{-3})$, to the fact that Weyl quantification
  $\Op(-if(x,\xi)$) is  bounded and to the definition of $\chi_1$, the operator $$K:={\widehat H}_{\te} -H_{\te}=\chi_1\left(\chi(-h^2\De+x^2)\Op(-if(x,\xi))\right)\chi_1$$ is of finite rank $O(h^{-3})$ and compactly supported. 
\hfill{$\square$}
\subsection{Upper bound for the number of resonances}\label{submajoration}
In this section we establish the upper bound of the number of resonances given in Theorem \ref{thmmajoration}.
\begin{lem}\label{lemmaj} Let $\rho>0$, $\Om$ be an open complex relatively compact subset of $\C$ and $H_{\te}$ be define as above. There exists $g$ satifying ({$\bf A_g$}) such that for $h$ small enough and $z\in {\overline \Om}\cap\{\emph{Im}\,z\geq \rho>0\}$,  we have $(z-H_{\te})^{-1}=O(1).$
\end{lem}
\Proof
We again use the notations of Section \ref{secdistortion}: $\phi_{\te}(x)=x+\te g(x)$ with $g(x)=0$ in the compact $B(0,R_0)$, and the notations of Subsection \ref{Appendix A} concerning $h_{\te}(x,\xi)$, 
 $U,\,U^{-1},\;d_{\te}$ and ${\Ti V}_{\te}(x,\xi)$ which satisfy
 \[U^{-1}h_{\te}(x,\xi)U=d_{\te}+{\Ti V}_{\te}(x,\xi).\]
The matrice $h_{\te}(x,\xi)$ is the semi-classical principal symbol of $H_{\te}$.

According to Section \ref{secdefinition}, the resonances are independent of the family  $U_{\theta}$. Then we can assume that $g(x)=0$ in the ball $B(0,R_g)\supset B(0,R_0)$, with $R_g>0$, sufficiently large such that 
$$\forall x\in\R^3,\;\;|x|>R_g>0,\;\;\|{\Ti V}_{\te}(x,\xi)\|\leq\frac{\rho}{2}.$$
Repeating  arguments of Subsection \ref{Appendix A}, we can prove that $(d_{\te}+{\Ti V}_{\te}(x,\xi)-z)$ is invertible, then $(h_{\te}(x,\xi)-z)$ is invertible and 
\begin{eqnarray} \nonumber
\|\p_x^\al \p_{\xi}^{\be}\left(h_{\te}(x,\xi)-z\right)^{-1}\|< C\<\xi\rangle^{-1-\be}.
\end{eqnarray}
Since we have: $$\|\p_x^\al \p_{\xi}^{\be}\left(h_{\te}(x,\xi)-z\right)\|< C\<\xi\rangle^{+1-\be},$$
the composition theorem of $h$-pseudo-differential operators implies 
 $$\Op(\left(h_{\te}(x,\xi)-z\right)^{-1})\Op(h_{\te}(x,\xi)-z)=1+O(h),$$
 where $O(h)$ corresponds to the norm in ${\cal L}(L^2)$.
\hfill{$\square$}\\\\
{\bf Proof of the Theorem \ref{thmmajoration}.}\\
Let ${\widehat K}(z)=K(z-{\widehat H}_{\te})^{-1}$ with ${\widehat H}_{\te}$, $K$ defined in Proposition \ref{propconst}. We remark that\\
$$(I+{\widehat K}(z))(z-{\widehat H}_{\te})=(z-{\widehat H}_{\te})+K=z-H_{\te}.$$
Thus, 
 the resonances $z\in$ Res($H)\cap \Om$ repeated with their multiplicities coincide with the zeros of the function
$$D(z)=\mbox{det}(I+{\widehat K}(z)).$$
%
Since $K$ is bounded and is of finite rank $O(h^{-3})$, we have
\[|D(z)|\leq e^{\|{\widehat K}(z)\|_{tr}}\leq e^{C_0h^{-3}}, \;\mbox{for all}\; z\in {\overline \Om}.\]
Using Lemma \ref{lemmaj}, we get $(z-H_{\te})^{-1}=O(1)$ for $\mbox{Im}  z\geq \rho>0$ and  $z\in {\overline \Om}$. Since
\begin{eqnarray}\label{KH}
(I+{\widehat K}(z))^{-1}=(z-{\widehat H}_{\te})(z-H_{\te})^{-1},
\end{eqnarray}
then
$$\|(I+{\widehat K}(z))^{-1}\|\leq C_1,\;\;\;\mbox{Im}  z\geq \rho>0.$$
 Writing the operator  $(I+{\widehat K}(z))^{-1}$ in the form $$(I+{\widehat K}(z))^{-1}=I-{\widehat K}(z)(I+{\widehat K}(z))^{-1},$$ we obtain
\[|\mbox{det}\left((I+{\widehat K}(z))^{-1}\right)|\leq e^{C_2h^{-3}},\;\; \mbox{Im}\,  z \geq \rho,\] 
 which implies  \[|D(z)|\geq Ce^{-C_3h^{-3}},z\in{\overline \Om}\cap \{\mbox{Im}\,  z\geq\rho\}.\] 
Now, applying Jensen's inequality in a slightly larger domain, we obtain Theorem $\ref{thmmajoration}$.\hfill{$\square$} 

\section{Representation of derivative of the spectral shift function }\label{secrep}
In this section we prove our principal result given in Theorem \ref{thmxi} and a generalization (see Theorem \ref{thmxig}). Moreover, we give a Sj\"ostrand type local trace formula.

The spectral shift function  $\xi(\la,h)\;(\in {\cal D}'(\R))$  associated to $H_0,H_1$ is defined
 (see \cite{VBDR}, \cite{DR3}, \cite{DY}) by \[\<\xi'(\la,h),f(\la)\rangle=\tr(f(H_1)-f(H_0)), \;\;\; f\in C_0^{\infty}(\R).\]
 In the following, we will use the notations:
  \[H_1=H,\;\;\;K_1\!:=\!K\!=\!{\widehat{H}}_{1,\te}-{H}_{1,\te}\!:=\!{\widehat{H}}_{\te}-H_{\te}\;\mbox{and}\;\;[a_\.\,]_0^1=a_1-a_0.\]
For an integer $m>3$,  we define the functions:
\begin{eqnarray}\label{eqsig}
\si_{\pm}(z)=(z^2+1)^m \tr\Big[(H_\.-i)^{-m}(H_\.+i)^{-m}(z-H_\.)^{-1}\Big]_0^1,\;\;\;\pm \mbox{Im}\,  z>0.
\end{eqnarray}
The $\si_{\pm}$ satisfy the relation
\begin{eqnarray}
\si_-(z)=\overline{\si_+(\bar z)},\;\;\;\mbox{Im}\,  (z)<0.
\end{eqnarray}
\begin{prop}For a potential $V$ satisfying the assumption $({\bf A_V})$ with $\de>3,$
the function $\te\longmapsto\Big[(H_{\.,\te}-i)^{-m}(H_{\.,\te}+i)^{-m}(z-H_{\.,\te})^{-1}\Big]_0^1$ is holomorphic from $D_{\epsilon}^+$ to the space of trace class operators. Moreover, for any $\te\in D_{\epsilon}^+$, we have
 \begin{eqnarray}\label{eqsigma}
\si_{\pm}(z)=
(z^2+1)^m \emph{tr}\Big[(H_{\.,\te}-i)^{-m}(H_{\.,\te}+i)^{-m}(z-H_{\.,\te})^{-1}\Big]_0^1,\;\;\;\pm \mbox{Im}\,  z>0.\;
\end{eqnarray} 
\end{prop}
\Proof
For $\te\in\R$, the operator $$\;(H_\.-i)^{-m}(H_\.+i)^{-m}(z-H_\.)^{-1},$$  
is unitarly equivalent to the operator 
$$(H_{\.,\te}-i)^{-m}(H_{\.,\te}+i)^{-m}(z-H_{\.,\te})^{-1}.$$
 Using the cyclicity of the trace, we deduce 
\begin{eqnarray}
\;\;\;\si_{\pm}(z)=(z^2+1)^m \tr\Big[(H_{\.,\te}-i)^{-m}(H_{\.,\te}+i)^{-m}(z-H_{\.,\te})^{-1}\Big]_0^1,\;\;\;\pm \mbox{Im}\,z>0,\;\te\in\R.
\end{eqnarray}
According to the proof of Theorem $\ref{thmdef}$,  the resolvent $(z-H_{\.,\te})^{-1}$ is analytic for $\te\in D_{\epsilon}^+$ and $\;z\in\Om\cap\{\mbox{Im}  z>0\}$.
Then, the function $\te\longmapsto(H_{\.,\te}-i)^{-m}(H_{\.,\te}+i)^{-m}(z-H_{\.,\te})^{-1}$ 
is also analytic on  $ D_{\epsilon}^+$.\\
Now, we treat the difference 
\begin{eqnarray}\label{eqABC}
\Big[(z-H_{\.,\te})^{-1}(H_{\.,\te}-i)^{-m}(H_{\.,\te}+i)^{-m}\Big]_0^1&=&A_1B_1C_1-A_0B_0C_0\nonumber\\
&=&A_1B_1(C_1-C_0)+A_1(B_1-B_0)C_0\nonumber\\&+&(A_1-A_0)B_0C_0.
\end{eqnarray}
Clearly, the terms $A_\.:=(z-{H}_{\.,\te})^{-1}$ for Im $z>0$, $B_\.:=(H_{\.,\te}-i)^{-m}$ and $C_\.:=(H_{\.,\te}+i)^{-m}$  are bounded.

For any integer $m>3,$ the term 
\begin{eqnarray}\label{eqdifff}
B_1(C_1-C_0)=\left(B_1(C_1-C_0)\<x\rangle^{\de}\<h\nabla_x\rangle^{m}\right)\left(\<h\nabla_x\rangle^{-m}\<x\rangle^{-\de}\right),
 \end{eqnarray}
 is analytic for $\te\in D_{\epsilon}^+$ with values in the space of trace class operators. This can be proved
using functional calculus in the framework of $h$-pseudo-differential operators (see \cite{MDJS}): The first factor $B_1(C_1-C_0)\<x\rangle^{\de}\<h\nabla_x\rangle^{m}$ is analytic for $\te\in D_{\epsilon}^+$, the second factor $\left(\<h\nabla_x\rangle^{-m}\<x\rangle^{-\de}\right)$ is 
in the space of trace class operators and its trace norm is bounded by $O(h^{-3})$. Then, 
the left-hand side 
of the equation (\ref{eqdifff}) is in the space of trace class operators and
  its trace norm is bounded by $O(h^{-3})$. The same argument can be used for the terms $A_1(B_1-B_0)$ and $(A_1-A_0)B_0$, then their trace norm are bounded by $O(h^{-3})$.

Since the function $\tr\Big[(H_{\.,\te}-i)^{-m}(H_{\.,\te}+i)^{-m}(z-H_{\.,\te})^{-1}\Big]_0^1$ is analytic with respect to $\te\in D_{\epsilon}^+$ and independent of $\te$ on the real axis, 
 the formula (\ref{eqsigma}) follows.\hfill{$\square$}
\\

Repeating the construction of ${\widehat{H}}_{1,\te}$, we can construct an operator ${\widehat{H}}_{0,\te}:D(H_0)\to \HH$ with the properties of ${\widehat{H}}_{0,\te}$ such that
the difference $K_0:={\widehat H}_{0,\te} -H_{0,\te}$ satisfies the properties of $K_1$ (see Proposition \ref{propconst}).
\begin{prop}\label{propsigma}
There exists a function $a_{+}(z,h)$  holomorphic in  $\Om$, such that for all \\$z\in\Om\cap\{\emph{Im} (z)>0\}$, we have:
\begin{eqnarray}
\si_{+}(z)=\emph{tr}\Big[(H_{\.,\te}-z)^{-1}K_\.({\widehat H}_{\.,\te}-z)^{-1}\Big]_0^1+a_{+}(z,h),
\end{eqnarray}
\[|a_+(z,h)|\leq C(\Om)h^{-3},\;\;\;z\in\Om,\]
with $C(\Om)$ a constant independent of $h$.
\end{prop}
\Proof
For $z\in\Om\cap\{\mbox{Im} z>0\}$, we have 
\begin{eqnarray}\label{eqdif}
(H_{\.,\te}-z)^{-1}=({\widehat{H}}_{\.,\te}-z)^{-1}+(H_{\.,\te}-z)^{-1}K_\.({\widehat{H}}_{\.,\te}-z)^{-1}.
\end{eqnarray}
From the equations $(\ref{eqdif})$ and $(\ref{eqsigma})$, we deduce:\\
$\si_+(z)=\left((z-i)(z+i)\right)^m \tr\Big[\left((\widehat{H}_{\.,\te}-z)^{-1}(H_{\.,\te}-i)^{-m}(H_{\.,\te}+i)^{-m}\right)\Big]_0^1$\\ \\
$\;\;\;\;\;\;\;\;\;\;\;\;\;\;\;+\left((z-i)(z+i)\right)^{m}\tr\Big[\left((H_{\.,\te}-z)^{-1}K_\.(\widehat{H}_{\.,\te}-z)^{-1}(H_{\.,\te}-i)^{-m}(H_{\.,\te}+i)^{-m}\right)\Big]_0^1$\\ \\
 $\;\;\;\;\;\;\;\;\;\;=A(z)+B(z).$\\
Starting with the resolvent equation, we obtain:
\begin{eqnarray}
\left((z-i)(z+i)\right)^m\hspace{-5mm}&\!\!\!\!\!\!\!\!\!\!\!&\hspace{-5mm}(H_{\.,\te}-i)^{-m}(H_{\.,\te}+i)^{-m}(H_{\.,\te}-z)^{-1}\nonumber\\\nonumber&=&(H_{\.,\te}-z)^{-1}
-\sum_{k=1}^m(z+i)^{k-1}(H_{\.,\te}+i)^{-k}\\\nonumber&-&(z+i)\sum_{k=1}^m(z-i)^{k-1}(H_{\.,\te}+i)^{-m}(H_{\.,\te}-i)^{-k}.
\end{eqnarray}
Using the last equation, the cyclicity of the trace and Propositon \ref{propconst}  
we obtain
\[B(z)=\tr\Big[K_\.(\widehat{H}_{\.,\te}-z)^{-1}\Big((H_{\.,\te}-z)^{-1}-\sum_{k=1}^m(z+i)^{k-1}(H_{\.,\te}+i)^{-k}\]\[-(z+i)\sum_{k=1}^m(z-i)^{k-1}(H_{\.,\te}+i)^{-m}(H_{\.,\te}-i)^{-k}\Big)\Big]_0^1\]
\[=\tr\Big[(H_{\.,\te}-z)^{-1}K_\.({\widehat H}_{\.,\te}-z)^{-1}\Big)\Big]_0^1+b(z).\;\;\;\;\;\;\;\;\;\;\;\;\;\;\;\;\;\;\;\;\] 
Since the operator $(\widehat{H}_{\.,\te}-z)^{-1}$ is bounded  and holomorphic in  $\Om$ by construction,  $b(z)$ is holomorphic and bounded by $O(h^{-3})$.

It remains to show that  \begin{eqnarray*}A(z)&=&\left((z-i)(z+i)\right)^m \tr\Big[(\widehat{H}_{\.,\te}-z)^{-1}(H_{\.,\te}-i)^{-m}(H_{\.,\te}+i)^{-m}\Big]_0^1\\
&=&\left((z-i)(z+i)\right)^m\tr(\widehat{A}_1B_1C_1-\widehat{A}_0B_0C_0),\end{eqnarray*} is holomorphic and bounded by $O(h^{-3})$.

We recall that the terms $\widehat{A}_\.:=(\widehat{H}_{\.,\te}-z)^{-1}$ for $z\in\Om$, $B_\.:=(H_{\.,\te}-i)^{-m}$ and $C_\.:=(H_{\.,\te}+i)^{-m}$  are bounded.
Using the assumption ({$\bf A_V$}) with $\de>3$, we treat the difference 
$(\widehat{A}_1B_1C_1-\widehat{A}_0B_0C_0)$ 
as (\ref{eqABC}). The only difference is for the term $(\widehat{A}_1-\widehat{A}_0)B_0$.
We write
\begin{eqnarray*}(\widehat{A}_1-\widehat{A}_0)B_0&=&(\widehat{H}_{1,\te}-z)^{-1}(\widehat{H}_{0,\te}-\widehat{H}_{1,\te})(\widehat{H}_{0,\te}-z)^{-1}(H_{0,\te}-i)^{-m},\\ 
\mbox{with }\;\;\widehat{H}_{0,\te}-\widehat{H}_{1,\te}&=&{H}_{0,\te}-{H}_{1,\te}+K_0-K_1.\end{eqnarray*}
Then, modulo a trace class operator uniformly bounded, with trace norm bounded by $O(h^{-3})$, we have
\begin{eqnarray*}(\widehat{A}_1-\widehat{A}_0)B_0&=&(\widehat{H}_{1,\te}-z)^{-1}\circ\left(({H}_{0,\te}-{H}_{1,\te})(H_{1,\te}-i)^{-m}\right)\\
&&\circ\left((H_{1,\te}-i)^{m}(\widehat{H}_{0,\te}-z)^{-1}(H_{0,\te}-i)^{-m}\right).\end{eqnarray*}
 The second factor $({H}_{0,\te}-{H}_{1,\te})(H_{1,\te}-i)^{-m}$ is trace class and its trace is $O(h^{-3})$, 
the first and the third factor are bounded.
Then, the term $(\widehat{A}_1-\widehat{A}_0)B_0$ is analytic for $z\in \Om$ with values in the space of trace class operators and its trace is bounded by $O(h^{-3})$ 
 and so is the difference $(\widehat{A}_1B_1C_1-\widehat{A}_0B_0C_0)$.
\hfill{$\square$}

\begin{lem}\label{lem}
For $f\in C_0^{\infty}(\R)$, we have
\begin{eqnarray}\label{eqxi,f}
\<\xi',f\rangle=\lim_{\ep\to0}\frac{i}{2\pi}\int f(\la)[\si_+(\la+i\ep)-\si_-(\la-i\ep)]d\la.
\end{eqnarray}
 This limit is taken in the sense of distribution.
\end{lem}
\Proof
We follow the proof of \cite[Lemma 1]{MDVP}. Let $f\in C_0^{\infty}(\R)$, $\ti f(z)\in C_0^{\infty}(\R^2)$ 
be an almost analytic extension of $f$ and \[g(x)=f(x)(x^2+1)^m.\] 
Then
\[g(H_{\.})=-\frac{1}{\pi}\int \bar\p_z\ti f(z)(z^2+1)^m(z-H_{\.})^{-1}L(dz),\]
where $L(dz)$ is the Lebesgue measure on $\C$. Clearly
\begin{eqnarray*}f(H_{\.})&=&(H_{\.}-i)^{-m}(H_{\.}+i)^{-m}g(H_{\.}) \\
&=&-\frac{1}{\pi}\int \bar\p_z\ti f(z)(z^2+1)^m(H_{\.}-i)^{-m}(H_{\.}+i)^{-m}(z-H_{\.})^{-1}L(dz)\end{eqnarray*}  which implies:
\begin{eqnarray}\label{eqtr}\tr\left(f(H_{1})-f(H_{0})\right)&=&-\frac{1}{\pi}\int \bar\p_z\ti f(z)(z^2+1)^m\nonumber\\
&\times& \tr\Big[(H_{\.}-i)^{-m}(H_{\.}+i)^{-m}(z-H_{\.})^{-1}\Big]_0^1L(dz).\;\;\;\;\;\;\;
\end{eqnarray}
 We have  $\si_{\pm}(z)=O(h^{-3}|\mbox{Im}  \,z|^{-2})$ and the derivative $\bar \p_z\ti f=O(|\mbox{Im} \,z|^N)\;$ for all $N\in\N$ ($f\in C_0^{\infty}(\R)$),  so we write the right-hand side of  $(\ref{eqtr})$ as
\[\<\xi',f\rangle=\tr\left(f(H_{1})-f(H_{0})\right)\]
\[=-\frac{1}{\pi}\lim_{\ep\to0}\left(\int_{\mbox{Im} \,z>0} \bar\p_z\ti f(z)\si_+(z+i\ep)L(dz)+\int_{\mbox{Im} \,z<0} \bar\p_z\ti f(z)\si_-(z-i\ep)L(dz)\right).\]
According to  Proposition $\ref{propsigma}$, the functions  $\si_+(z+i\ep)$ and $\si_-(z-i\ep)$ are holomorphic  in  $\{z\in \Om;\;\mbox{Im}  z>0\}$ and $\{z\in \Om;\;\mbox{Im}  z<0\}$ respectively.  Applying the Green formula, we obtain the lemma.\hfill{$\square$}\\ \\
 Before the proof of Theorem  $\ref{thmxi}$, let us give the following proposition:
\begin{prop}\label{propomegatilde}$($see \cite{JS1}, \cite{JS2}$)$ Let
 $F(z,h)$ be a holomorphic function in an open simply connected domain $\Om$ containing  
 a number $N(h)$ of zeros.
We suppose that,\[F(z,h)=O(1)e^{O(1)N(h)},\;\;z\in\Om,\]
and for all $\rho>0$  small enough,  there exists  $C>0$ such that for all $z\in\Om_{\rho}:=\Om\cap\{\emph{Im} \,z>\rho
\}$ we have \[|F(z,h)|\geq e^{-CN(h)}.\]
 Then for all open simply connected subset  $\ti\Om\Subset\Om$ there exists $g(.,h)$ holomorphic in $\ti\Om$ such that \[F(z,h)=\prod_{j=1}^{N(h)}(z-z_j)e^{g(z,h)},\;\;\p_zg(z,h)=O(N(h)),\;\;z\in\ti\Om.\]

\end{prop}
{\bf Proof of Theorem \ref{thmxi}.}
 We follow the argument of  Sj\"ostrand (\cite{JS2}). Let $${\widehat K}_\.(z)=K_\.(z-{\widehat H}_{\.,\te})^{-1}.$$
As  \cite[Equation(4.31)]{JS2} we have the representation
\begin{eqnarray*}-\tr\left((H_{\.,\te}-z)^{-1}K_\.({\widehat H}_{\.,\te}-z)^{-1}\right)&=&\tr\left(\left(1+\widehat{K}_\.(z)\right)^{-1}\p_z\widehat{K}_\.(z)\right)\\ \;\;\;\;\;\;\;\;\;\;\;&=&\p_z\mbox{log det}\left(1+\widehat{K}_\.(z)\right).\end{eqnarray*}
 From Subsection  $\ref{submajoration}$  the resonances are the zeros of the function 
\[D(z,h)=\mbox{det}\left(I+{\widehat K_1}(z)\right)=O(1)e^{ch^{-3}}.\]
Since the function $\mbox{det}(1+\widehat{K}_0(z))$ has no zeros in $\Om$ (see (\ref{KH}) and Remark \ref{remark6})
, the term $\p_z\mbox{log det}\left(1+\widehat{K}_0(z)\right)$ is analytic and using the Proposition \ref{propomegatilde}, it is bounded by $O(h^{-3})$.

We recall the $\Res(H)$ be the set of resonances of $H$ and let \[D(z,h)=G(z,h)\prod_{w\in \Res(H)\cap\Om}(z-w),\]
where, $G(z,h)$  and its inverse are holomorphic functions in  $\Om.$  Obviously,
\begin{eqnarray}\label{eqsome}
\p_z\mbox{log}\,D(z,h)=\p_z\mbox{log}\,G(z,h)+\sum_{w\in \Res(H)\cap\Om}\frac{1}{z-w}. 
\end{eqnarray}
 Using  Proposition $\ref{propomegatilde}$, we have\[ \;\;\;\;\left|\p_z\mbox{log}\,G(z,h)\right|\leq C(\ti\Om)h^{-3},\;\;\;z\in\ti\Om,\]
where $\ti\Om\subset\subset\Om$  is an open simply connected set and  $C(\ti\Om)$ is independent of  $h$.\\ \\
 Now, we treat the non-holomorphic term in  $\left(\si_+(\la+i\ep)-\si_-(\la-i\ep)\right)$  when $\ep\to 0,$ which is
\[\sum_{w\in \Res(H)\cap\Om}\left(\frac{1}{\la+i\ep-w}-\frac{1}{\la-i\ep-\overline{w}}\right), \;\;\mbox{for} \la\in I.\]
If $\mbox{Im} (w)\neq0,$ we have \[\frac{-1}{2i\pi}\lim_{\ep\to0}\left(\frac{1}{\la+i\ep-w}-\frac{1}{\la-i\ep-\overline{w}}\right)=\frac{-\mbox{Im} (w)}{\pi|\la-w|^2},\]
while for $w\in\R$ we get\[\hspace{13mm}\frac{-1}{2i\pi}\lim_{\ep\to0}\left(\frac{1}{\la+i\ep-w}-\frac{1}{\la-i\ep-{w}}\right)=\de(\la-w)=\de_w(\la).\]
The second limit is taken in  the sense of distributions. \\
Lemma $\ref{lem}$ and Proposition $\ref{propsigma}$ 
show that the function  $r(z,h)=g(z,h)-{\overline g(\bar z,h)}$, with $g(z,h)=a_+(z,h)+\p_z\mbox{log}\,G(z,h)+\p_z\mbox{log det}\left(1+\widehat{K}_0(z)\right)$ is a holomorphic function in $\Om$ and satisfies the following estimate:
\begin{eqnarray}
|g(z,h)|\leq C(\Om)h^{-3}, \;\;z\in W,
\end{eqnarray}
with $C(\Om)>0$ independent of $h$.\hfill{$\square$}\\

Theorem \ref{thmxi} can be extended to a more general situation:
\begin{thm}\label{thmxig} Assume that $H_1=H_0+V_1,\;H_2=H_0+V_2$. The potentials  $V_1,V_2$ (resp. $V=V_1-V_2$) satisfy the assumption ${\bf (A_V)}$ with $\de>0$ (resp. $\de>3$). Let $\Om$ be a complex domain satisfying the assumption ${\bf(A_{\Om}^{\pm})}$, $W\Subset \Om$ be an open simply connected and relatively compact set which is symmetric with respect to $\R$. Assume that $I=W\cap \R$ is an interval. Then for all $\la\in I$  we have a representation of the derivative of the spectral shift function associated to the pair $(H_2,\ H_1)$ of the form:
\begin{eqnarray}\label{eqderiveg}
\xi'(\la,h)=\frac{1}{\pi}\emph{Im}  \,r(\la,h)+{\Big[}\sum_{\begin{array}{c}
\scriptstyle w\in \res(H_\.)\cap\Om \\
\scriptstyle \emph{Im}\,  w\neq0
\end{array}}\frac{-\emph{Im}\,  w}{\pi|\la-w|^2}+\sum_{w\in \res(H_\.)\cap I}\de_w(\la){\Big]_2^1},
\end{eqnarray}
where $r(z,h)=g(z,h)-\bar g( \bar z,h)$, $g(z,h)$ is a holomorphic function in  $\Om$ which satisfies the following estimate: 
\begin{eqnarray}
|g(z,h)|\leq C(W)h^{-3}, \;\;z\in W,
\end{eqnarray}
with $C(W)>0$  independent of  $h$. Here $\de_w(\cdot)$ is the Dirac 
mass at $w \in \R$. 
\end{thm}
{\bf 
 Proof.} We denote $H_{2,\te}=U_{\te}H_2U_{\te}^{-1}$ ($U_{\te}$ defined in Section $\ref{secdistortion}$). As in Subsection \ref{Appendix A}, one constructs ${\widehat H}_{2,\te}:D(H)\longmapsto \HH$  with the following properties:\\ $K_2:={\widehat H}_{2,\te} -H_{2,\te}$  is of finite rank  $O(h^{-3})$, has a compact support in the sense that $K_2=\chi_2 K_2\chi_2$ if $\chi_2 \in \C_0^{\infty}$ is equal to 1 on $B(0,R)$ for some sufficiently large $R$, and \[({\widehat H}_{2,\te}-z)^{-1}=O(1):\HH \longmapsto D(H),\; \mbox{uniformly for}\; z\in {\overline \Om}.\]
We repeat the proof of the Theorem $\ref{thmxi}$ replacing $K_0$ by $K_2$ and ${\widehat K}_0(z)$ by ${\widehat K}_2(z)=K_2(z-{\widehat H}_{2,\te})^{-1}$. Consequently $\p_z\mbox{log det}\left(1+\widehat{K}_0(z)\right)$ is replaced by $\p_z\mbox{log det}\left(1+\widehat{K}_2(z)\right)$ which is a non-holomorphic function. We treat this term as the term $\p_z\mbox{log det}\left(1+\widehat{K}_1(z)\right)$ in the proof of Theorem $\ref{thmxi}$.\hfill{$\square$}
 \begin{rem}
 The equation $(\ref{eqderive})$ shows that the spectral shift function $\xi(\la,h)$ satisfies
\begin{eqnarray} \label{eqxi}
\xi(\la,h)-\xi(\la_0,h)&=&\sum_{\begin{array}{c}
\scriptstyle w\in \res(H_1)\cap\Om \\
\scriptstyle \emph{Im}\,  w\neq0
\end{array}}\frac{1}{\pi}\int_{\la_0}^{\la}\frac{-\emph{Im}\,  w}{|\mu-w|^2}d\mu+\frac{1}{\pi}\int_{\la_0}^{\la}\emph{Im} \; r(\mu,h)d\mu\nonumber\\ 
&+&\;\,\#\{\mu\in[\la_0,\la];\;\mu\in\si_{d}(H_{1})\}
.\end{eqnarray}
In particular, for  $\la \in I\backslash\si_d(H_{1})$ the distribution $\xi(\la,h)$  is continuous, and the function \[\eta(\la,h)-\eta(\la_0,h)=\xi(\la,h)-\xi(\la_0,h)-\,\#\{\mu\in[\la_0,\la];\;\mu\in\si_{d}(H_{1})\}\] is real analytic in  $I$.\hfill{$\square$} 
\end{rem}

Repeating the argument used in the proof of \cite[Theorem 4]{VBVP}, 
 the following theorem is a direct consequence of Theorem  $\ref{thmxi}$.
\begin{thm}\label{thmtrace}{\bf(Local trace formula)}
Let $\Om$ be an open complex, simply connected and relatively compact set satisfying the assumption ${\bf (A_{\Om}^{\pm})}$ such that $I=\Om\cap\R$ is an interval.\\
 We suppose that $f$ is a holomorphic function in  $\overline{\Om}$ and $\psi\in C_0^{\infty}(\R)$ satisfies 
\[\psi(\la)=\left\{\begin{array}{rl}
0,&\;\;d(I,\la)>2\ep,\\
1,&\;\;d(I,\la)<\ep,\end{array}\right.\]   
  where $\ep>0$
and sufficiently small. Then 
  \[tr[(\psi f)(H_\.)]_0^1=\sum_{z\in \res(H_1)\cap\Om}f(z)+E_{\Om,f,\psi}(h), \mbox{  with}\]
\[|E_{\Om,f,\psi}(h)|\leq M(\psi,\Om)\sup\{|f(z)|;\;0\leq d(\p\Om,z)\leq 2\ep,\; \emph{Im}\, (z)\leq0\}h^{-3}.\]
\end{thm}
\section{Weyl Asymptotics }\label{secweyl}
In this section 
we obtain a Weyl-type asymptotic for the spectral shift function $\xi(\la,h)$ associated to operators $H_0$ and $H_1=H_0+V$.  Here 
we assume that $V$ is an electro-magnetic potential ($\ref{eqV}$),
\[H_1=-\sum_{j=1}^3\al_j(ich\p_j+eA_j)+mc^2\be+ev.\]
In the following, we fix $I_0\subset\R\backslash\{\pm mc^2\}$ and choose $W_0$ an open simply connected, relatively compact subset of $\Om$ satisfying the assumption $({\bf A_{\Om}^+})$ such that $I_0=W_0\cap\R$.

For the $h$-pseudo-differential and functional calculus
 for the Dirac operator, we refer to (\cite{MDJS},\cite{VBDR},\cite{DR1},\cite{BHDR}). We recall that $H_{\nu}=\Op(\DD_{\nu})$ and $\vphi(H_\nu)$ are $h$-pseudo-differential operators for a smooth function $\vphi$. The semi-classical symbol $\DD_{\nu}$ is defined in $(\ref{eqDnu}).$\\

Let us introduce the intervals $I_1,\ I_2\subset I_0$  neighborhoods  of $\la_1,\;\la$ respectively such that, each $\la\in I_1\cup I_2$ is a noncritical energy level for $H$ (see Definition \ref{defnoncritical}).
Let $\vphi_j\in C_0^{\infty}(\R,\R^+) $ 
 be such that  
\begin{eqnarray}\label{eqpartition}
\vphi_1=1\;\mbox{on}\; I_1,\;\;\;\vphi_2=1\; \mbox{on}\; I_2\;\mbox{and} \;\;\;\vphi_1+\vphi_2+\vphi_3=1\; \mbox{on}\; I_0.
\end{eqnarray}

Consider a function $\te(t)\in C_0^{\infty}(]-\de_1,\de_1[)$, $\te(0)=1,$ $\te(-t)=\te(t)$, so that the Fourier transform  $\hat\te$ of $\te$ satisfies $\hat\te(\la)\geq0 \;$on$\;\R$, and assume that there exist  $0<\epsilon_0<1,\;\de_0>0$, such that $\hat\te(\la)\geq\de_0>0 \;\mbox{for}\;|\la|\leq\epsilon_0$. Next, we introduce
\[({\cal F}_h^{-1}\te)(\la)=(2\pi h)^{-1}\int e^{it\la h^{-1}}\te(t)dt=(2\pi h)^{-1}\hat\te(-h^{-1}\la).\] 
To prove Theorem \ref{thmweyl}, we need the proposition:
\begin{prop}\label{proptrace}
For the trace involving $H_{\nu},\;\;\nu=0,1,$ we have for $\la\in I_j,$
\begin{eqnarray}\label{eqw1}
\emph{tr}\left(\Big[({\cal F}_h^{-1}\te)(\la-H_\.)\vphi_j(H_\.)\Big]_0^1\right)
=w_j(\la)h^{-3}+O(h^{-2}),\;\;j=1,2, 
\end{eqnarray}
with $w_j(\la)\in C_0^\infty(I_j)$ and $O(h^{-2})$ uniform with respect to $\la\in I_j$.
\end{prop}
\Proof
Proposition \ref{proptrace} is closed to the calculation of the trace in \cite[Section 4]{JB} and to the appendix of \cite{VBVP} for the Schr\"odinger operator. But, here we use a trick of Robert \cite{VBDR}. We fix $j=2$ (it is similar for $j=1$). The proof of ($\ref{eqw1}$) is obtained following these two steps:

$\bullet$ First, 
we recall that $\la\in I_2$ and $\mbox{Supp}\,\te(t)\subset[-\de_1,\de_1]$. Let us write

\begin{eqnarray*}{\cal T}=\tr\Big[({\cal F}_h^{-1}\te)(\la-H_\.)\vphi_2(H_\.)\Big]_0^1&\!\!\!\!\!=\!\!\!\!\!&\tr\Big[\!\!\int \frac{\te(t)}{2\pi h}e^{it(\la-H_\.)h^{-1} }\vphi_2(H_\.)dt\Big]_0^1\\
&\!\!\!\!\!=\!\!\!\!\!&\frac{1}{2\pi h}\int e^{it\la h^{-1}}\te(t) \tr[e^{-itH_\.h^{-1}}\vphi_2(H_\.)]_0^1dt
\end{eqnarray*} 
In the order to calculate the trace $$\tr(f(H_1)-f(H_0)),\;\;\;\mbox{for all}\;f\in C_0^\infty(\R\backslash\{\pm mc^2\})$$
we use \cite[Proposition 3.2]{VBDR}. If we note $W(h)=Q-\frac12i[Q,{\cal A}(h)]$ with $Q=H_1^2-H_0^2$, ${\cal A}(h)=\frac12(x\cdot h\p_x+h\p_x\cdot x)$ and $[Q,{\cal A}(h)]=Q{\cal A}(h)-{\cal A}(h)Q$, we have
\begin{eqnarray}\label{formulae}\tr(f(H_1)-f(H_0))=\tr(W(h)(H_1^2-m^2c^4)^{-1}f(H_1)).
\end{eqnarray}

Applying the formula (\ref{formulae})  for $f(\la)=e^{-it\la h^{-1}}\vphi_2(\la)$, we have
$${\cal T}=\frac{1}{2\pi h}\int e^{it\la h^{-1}}\te(t) \tr\left(W(h)(H_1^2-m^2c^4)^{-1}e^{-itH_1h^{-1}}\vphi_2(H_1)\right)dt.$$
\begin{rem}
 Of course $(H_1^2-m^2c^4)^{-1}$ is not well defined, however for $f\in C_0^\infty(\R\backslash\{\pm mc^2\})$, we can define $(H_1^2-m^2c^4)^{-1}f(H_1)$ as the self-adjoint operator $\vphi(H_1)$ where $\vphi\in C_0^\infty(\R)$ satisfies:
$$\vphi(\la)=\left\{\begin{array}{rl}
(\la^2-m^2c^4)^{-1}f(\la)&\emph{for  } \la\neq\pm mc^2,\\
0&\emph{for  } \la=\pm mc^2.
\end{array}\right.$$
\end{rem}

$\bullet$ Now, we treat ${\cal T}$  following the analysis of \cite[Section 4.2]{JB}.  
By $h$-pseudo-differential calculus, we obtain the existence of a $h$-pseudo-differential operator $S$ which is trace class with symbol\begin{eqnarray}\label{eqS}
s(x,y,\xi,h)\in {\cal S}^0(\<x\rangle^{-\de}\<\xi\rangle^{-N}),\;\;\;\forall N\in \N,\;\de>3,
\end{eqnarray}
having compact support in $\xi$ and in ($x-y$) (i.e. $\mbox{supp}_{(x-y)}(s)=\{x-y,\;\exists\; \xi\;;\,(x,y,\xi,h)\!\in \!\mbox{supp}(s)\}$ is compact) and support in $\{(x,\xi);\;|x|>R,\; (x,\xi)\in \DD_1^{-1}(I_2)\}$, with $\DD_{1}$ the semi-classical symbol of $H_{1}$, so that
\[
{\cal T}=\frac{1}{2\pi h}\tr\left(\int e^{it\la h^{-1}}\te(t)e^{-itH_1h^{-1}}Sdt\right)+O(h^{\infty}).\]

Using Theorem \ref{thmpropagator} in Appendix A and the hypothesis of $S$ by composition of Fourier integral operators,  we obtain a Fourier integral operator ${\Ti{\cal U}}_t={\Ti{\cal U}}_t^++{\Ti{\cal U}}_t^-$, such that for $|t|\leq \de_1$ and $\de_1$ sufficiently small, we have\begin{eqnarray}\label{eqnormtr}\|{\Ti{\cal U}}_t-e^{-itH_1h^{-1}}S\|_{tr}=O(h^{\infty}),\end{eqnarray}
where the kernel of the operator $\int e^{it\la h^{-1}}\te(t){\Ti{\cal U}}_t dt $ is equal to  ${\Ti K}^+(x,y;h)+{\Ti K}^-(x,y;h)$ with
\[{\Ti K}^{\pm}(x,y;h)=\frac{1}{(2\pi h)^3}\int\!\!\int e^{i(t\la+\Phi^{\pm}(t,x,\xi)-y\cdot\xi)h^{-1}}\te(t){\Ti E}^{\pm}(t,x,y,\xi;h)dtd\xi.\]
The amplitudes ${\Ti E}^{\pm},$ satisfy
 \[{\Ti E}^{\pm}(t,x,y,\xi;h)\in {\cal S}^0(\<x\rangle^{-\de}\<\xi\rangle^{-N}),\;\;\;\forall N\in \N\]
and is compactly supported in $\xi$ and in ($x-y)$.\\
Using the Taylor sum formula for the functions $\Phi^{\pm}(t,x,\xi)$  in a neighborhood of $t=0$, we have:
 \[\Phi^{\pm}(t,x,\xi)= x\cdot\xi-tH_1^{\pm}(x,\xi)+O(t^2).\]
  We will deduce that 
  ${\cal T}={\cal T}^++{\cal T}^-$, with  
\[
{\cal T}^{\pm}=\frac{1}{(2\pi h)^4}\int\!\!\int\!\!\int e^{i(t\la+\Phi^{\pm}(t,x,\xi)-x\cdot\xi)h^{-1}}\te(t){\Ti E}^{\pm}(t,x,x,\xi;h)dtdxd\xi+O(h^{\infty}).\]
Moreover, the symbol  ${\Ti E}^{\pm}(t,x,x,\xi;h)$ has support in $\{(x,\xi);\,|x|>R,\,|\xi|\leq C_1,\,(x,\xi)\in \DD_1^{-1}(I_2)\}$, so that for all $\al$ and $|t|\leq\de_1$, we have \begin{eqnarray}\label{eqA1}|\p^{\al}{\Ti E}^{\pm}(t,x,x,\xi;h)|\leq C_{\al}\<x\rangle^{-\de},\;\;\de>3.
\end{eqnarray}
The last estimate enables us to calculate 
${\cal T}$ by using an infinite partition of unity 
\[\sum_{\al\in\N^3}\Psi(x-\al)=1,\;\;\;\forall x\in \R^3,\]
where $\Psi\in C_0^{\infty}(K), \Psi\geq0,$ $K$ being a neighborhood of the unit cube. Consequently, for every fixed $h\in]0,h_0]$, we have 
\begin{eqnarray*}
{\cal T}^{\pm}&=&\frac{1}{(2\pi h)^4} \lim_{m\to\infty}\int\!\!\!\int\!\!\!\int\!e^{i(t\la+\Phi^{\pm}(t,x,\xi)-x\cdot\xi)h^{-1}}\te(t)\\&\times&\sum_{|\al|\leq m}\Psi(x-\al){\Ti E}^{\pm}(t,x,x,\xi;h)dtdxd\xi+O(h^{\infty})=\lim_{m\to\infty}I_m^{\pm}+O(h^{\infty}),
\end{eqnarray*}
and we reduce the problem to the analysis of the integrals $I_m^{\pm}$. Concerning the phase function, we observe that
\begin{eqnarray}\label{eqptcritics}
 t\la+\Phi^{\pm}(t,x,\xi)-x\cdot\xi=t(\la-H_1^{\pm}(x,\xi)+O(t)),
\end{eqnarray}
where $O(t)$ and its derivatives are uniformly bounded on the support of $\te(t){\Ti E}^{\pm}(t,x,x,\xi;h)$  since the derivatives of $(\Phi^{\pm}(t,x,\xi)-x\cdot\xi)$ are bounded on this set.

Now we look for critical points of the phase function $(t\la+\Phi^{\pm}(t,x,\xi)-x\cdot\xi)$. Putting the  derivative with respect to $t$  equal to $0$, we see that $H_1^{\pm}(x,\xi)=\la+O(t)$. Since $\p_{x,\xi}H_1^{\pm}(x,\xi)\neq0$, when $H_1^{\pm}(x,\xi)=\la$, and putting the derivative of the phase function $t(\la-H_1^{\pm}(x,\xi)+O(t))$ with respect to $ H_1^{\pm}(x,\xi)$ equal to $0$, we have
$$t=O(t^2).$$ Then the phase is critical for $|t|$ small precisely when   $t=0,\;\;\la=H_1^{\pm}.$ Near any such critical point we choose local coordinates $t,\,H_1^{\pm}(x,\xi),\,w_1,\cdots,w_5$ and consider the Hessian of (\ref{eqptcritics})
with respect to $t,\,H_1^{\pm}(x,\xi)$ at the critical point:
\begin{displaymath}
\left(\begin{array}{cc}
\star& -1 \\ 
-1 & 0
\end{array}\right)
\end{displaymath}
This is a non-degenerate matrix of determinant $-1$ and of signature $0$.
By the stationary phase method we obtain
\[I_m^{\pm}=\frac{\psi^{\pm}(\la)}{(2\pi h)^3}\int_{\la=H_1^{\pm}}\sum_{|\al|\leq m}\Psi(x-\al){\Ti E}^{\pm}(0,x,\xi,\la;h)L_{\la}^{\pm}(dw)+O(h^2),\]
where $L_{\la}^{\pm}(dw)$ is the Liouville measure on $\la=H_1^{\pm}$ and the remainder $O(h^{-2})$ is uniform with respect to $\la\in I_2$ and $m\in \N$. Here $\psi^{\pm}(\la)\in C_0^\infty(I_2)$. Taking the limit $\lim_{m\to\infty}I_m^{\pm}$, we obtain an asymptotics of 
${\cal T}.$\hfill{$\square$} 

\begin{lem} With the above definitions of  $\te(t),\;\xi(\la,h),\;\vphi_j(\la),\;I_j,\;\;j=1,2,$ we have
\begin{eqnarray}
\label{eqconv}\int_{-\infty}^{\la}{\cal F}_h^{-1}\te*\vphi_j\xi'(\mu,h)d\mu-\int_{-\infty}^{\la}\vphi_j(\mu)\xi'(\mu)d\mu=O(h^{-2}),\;\la\in I_j.
\end{eqnarray}
\end{lem}
\Proof
We deal only with the analysis of (\ref{eqconv}) for $j=2$ since that of $j=1$ is similar. According to Theorem $\ref{thmxi}$, there exists a holomorphic function $r(z,h)$ in $\Om$ such that for all $\la\in I_0=W_0\cap\R$,
  we have
\[\xi'(\la,h)=\frac{1}{\pi}\mbox{Im}  \,r(\la,h)+\sum_{\begin{array}{c}
\scriptstyle w\in \Res(H_1)\cap\Om \\
\scriptstyle \mbox{Im}  w\neq0
\end{array}}\frac{-\mbox{Im}  w}{\pi|\la-w|^2}+\sum_{w\in \Res(H_1)\cap I_0}\de_w(\la),\]
where $r(z,h)$  satisfies the following estimate:
\begin{eqnarray}\label{eqrrr}
|r(z,h)|\leq C(W)h^{-3}, \;\;z\in W,
\end{eqnarray}
with $C(W)>0$ independent of $h$. Let us denote
\[
G_{\vphi_2}(\la)=\frac{1}{\pi}\int_{-\infty}^{\la}\mbox{Im} \;r(\mu,h)\vphi_2(\mu)d\mu,\]
\begin{eqnarray}\label{eqm}
M_{\vphi_2}(\la)=\sum_{\begin{array}{c}
\scriptstyle w\in \mbox{Res}(H_1)\cap\Om \\
\scriptstyle \mbox{Im}  w\neq0
\end{array}}\int_{-\infty}^{\la}\frac{-\mbox{Im }  w}{\pi|\la-w|^2}\vphi_2(\mu)d\mu+\sum_{w\in \mbox{Res}(H_1)\cap ]c_0,\la]}\vphi_2(w).
\end{eqnarray}
Using the Cauchy inequalities and (\ref{eqrrr}), it  follows easily that \[G_{\vphi_2}'(\la)=O(h^{-3}) \;\mbox{and}\;G_{\vphi_2}^{\prime\prime}(\la)=O(h^{-3}),\]and we immediately obtain
\begin{eqnarray}\label{eqg}
{\cal F}_h^{-1}\te*G_{\vphi_2}'-G_{\vphi_2}'=O(h^{-2}).
\end{eqnarray}
Now, we want to apply a Tauberian theorem (see \cite[Theorem V-13 ]{DR1}) for the increasing function $M_{\vphi_2}(\la)$. For this purpose, we need the estimates
 \begin{eqnarray}\label{eqmphi}
 M_{\vphi_2}(\la)=O(h^{-3}),\;\;\;
 \frac{d}{d\la}({\cal F}_h^{-1}\te*M_{\vphi_2})(\la)=O(h^{-3}),\;\;\forall \la\in\R,
 \end{eqnarray}
 and the equality $\;\;M_{\vphi_2}(\mu)=G_{\vphi_2}(\mu)=0,\;\;\;\mu\leq inf I_2.$\\
 The first estimate in (\ref{eqmphi}) follows easily from the equation ($\ref{eqm}$) with the upper bound of the number of the resonances in $\Om$ (see Theorem \ref{thmmajoration}), and the second follows from ($\ref{eqw1}$) and the equation
  \[\frac{d}{d\la}({\cal F}_h^{-1}\te*M_{\vphi_2})(\la)={\cal F}_h^{-1}\te*\vphi_2\xi'(\la)-\frac{d}{d\la}({\cal F}_h^{-1}\te*G_{\vphi_2})(\la).\]
  Then, according to the Tauberian theorem we have 
  \[({\cal F}_h^{-1}\te*M_{\vphi_2})(\la)=M_{\vphi_2}(\la)+O(h^{-2}),\]
    this enabled us to obtain 
  \begin{eqnarray*}\label{eqmprime}
  \int_{-\infty}^{\la}\vphi_2(\mu)\xi'(\mu)d\mu&=&M_{\vphi_2}(\la)+\int_{-\infty}^{\la}G_{\vphi_2}'(\mu)d\mu\\
 &=&\int_{-\infty}^{\la} \frac{d}{d\mu}({\cal F}_h^{-1}\te*M_{\vphi_2}+{\cal F}_h^{-1}\te*G_{\vphi_2})(\mu)d\mu+O(h^{-2})\\
 &=&\int_{-\infty}^{\la}{\cal F}_h^{-1}\te*\vphi_2\xi'(\mu,h)d\mu+O(h^{-2}).\;
  \end{eqnarray*}
\hfill{$\square$}

{\bf Proof of Theorem \ref{thmweyl}.}
For $\la_1\in I_1,\;\la\in I_2$, using the functions defined in (\ref{eqpartition}),
we  have
\begin{eqnarray}\label{eqxiweyl}
\xi(\la,h)-\xi(\la_1,h)&=&\int_{-\infty}^{\la}\vphi_1(\mu)\xi'(\mu,h)d\mu-\int_{-\infty}^{\la_1}\vphi_2(\mu)\xi'(\mu,h)d\mu\nonumber\\ \nonumber
&-&\int_{-\infty}^{\la_1}\vphi_1(\mu)\xi'(\mu,h)d\mu+\int_{-\infty}^{\la}\vphi_2(\mu)\xi'(\mu,h)d\mu\\ 
&+&\int_{\la_1}^{\la}\vphi_3(\mu)\xi'(\mu,h)d\mu.
\end{eqnarray}

Since $\vphi_j=0\;$on$\;I_{3-j}$ for $j=1,2$, the first term (resp. the second term) is independent of $\la\in I_2$ (resp. $\la_1\in I_1$) and is equal to $\;\tr[\vphi_1(H_.)]_0^1=C(\vphi_1)h^{-3}+O(h^{-2})$ (resp. $\tr[\vphi_2(H_.)]_0^1=C(\vphi_2)h^{-3}+O(h^{-2})$), where $C(\vphi_j)$ is a constant depending on $\vphi_j$ for $j=1,2$. Since $\vphi_3=0$ on $I_j$, $j=1,2$, the last term is independent of $\la\in I_2$, $\la_1\in I_1$ and is equal to $C(\vphi_3)h^{-3}+O(h^{-2})$, where $C(\vphi_3)$ is a constant depending on $\vphi_3$. The proof of this results is based on the functional calculus in the framework of $h$-pseudo-differential operators. 
\\

Using the equations (\ref{eqw1}), (\ref{eqconv}) and  (\ref{eqxiweyl})
 we complete the proof of the asymptotic expansion (\ref{weyl}) by writing 
\begin{eqnarray}({\cal F}_h^{-1}\te*(\vphi_j\xi'))(\la)
&=&\,\<({\cal F}_h^{-1}\te)(\la-{\bf \cdot})\vphi_j({\bf \Large\cdot}),\xi'\rangle\nonumber\\
&=&\tr\left(\Big[({\cal F}_h^{-1}\te)(\la-H_\.)\vphi_j(H_\.)\Big]_0^1\right)\nonumber\\
&=&w_j(\la)h^{-3}+O(h^{-2}),\;\;j=1,2. 
\end{eqnarray}
 It remains to compute the Weyl term (\ref{eqpremierterm}).\\\\
According to the definition of the spectral shift function $\xi(\la,h)$ in $(\ref{defxi})$, we have:
\begin{eqnarray}\label{eq*}\<\xi'(\la,h),\vphi(\la)\rangle=\tr(\vphi(H_1)-\vphi(H_0)),\;\;\vphi(\la)\in C_0^{\infty}(\R).
\end{eqnarray}
We use weak asymptotics which is a direct consequence of functional calculus in the framework of  $h$-pseudo-differential operators, as settled in \cite{MDJS}, \cite{DR1}, \cite{VBDR}. We find
\[H_{\nu}=\Op(\DD_{\nu}),\;\;\;\nu=0,1,\;(\DD_{\nu}\;\mbox{defined in}\;(\ref{eqDnu})),\]
and
\begin{eqnarray*}\tr\left(\vphi(H_1)-\vphi(H_0)\right)&=& h^{-3}\sum_{j\geq0}\ga_j(\vphi)h^j\\
&=&h^{-3}\ga_0(\vphi)+O(h^{-2}),
\end{eqnarray*}
with $\ga_0(\vphi)=(2\pi)^{-3}\int_{\R^3}\!\!\int_{\R^3}\tr\left(\vphi(\DD_1(x,\xi))-\vphi(\DD_0(x,\xi))\right)dxd\xi$.\\ ($\tr(A)$ is the trace of the matrix $A$).\\\\
The matrix $\DD_{\nu}(x,\xi)$ is Hermitian and has two eigenvalues $H_{\nu}^{\pm}(x,\xi)$ (see ($\ref{eqH+-}$)), then
\[\tr\left(\vphi(\DD_1)-\vphi(\DD_0)\right)=2\left(\vphi(H_1^+)+\vphi(H_1^-)-\vphi(H_0^+)-\vphi(H_0^-)\right).\]
According to the asymptotic expansions (\ref{weyl}) 
 and (\ref{eq*})  we obtain $$w(\la,\la_1)=w(\la)-w(\la_1),$$ with
\[w(\la)=\frac{1}{4\pi^3}\!\int_{\R^3}\!\!\left(\!\int_{H_1^+(x,\xi)\leq\la}\!\!d\xi-\int_{H_0^+(x,\xi)\leq\la}\!\!d\xi-\int_{H_1^-(x,\xi)\geq\la}\!\!d\xi+\int_{H_0^-(x,\xi)\geq\la}\!\!d\xi\right)\!dx.\]
Putting $\zeta_{\nu}=c\xi-\nu eA(x)$ for $\nu=0,1$ and $\zeta=r\om\;(\om\in S^2)$, we get
\begin{eqnarray*}\pm H_{\nu}^{\pm}\leq\pm\la&\Leftrightarrow&\left(\zeta_{\nu}^2+(mc^2+\nu\frac{e(v_+-v_-)}{2})^2\right)^{\frac12}\pm\nu\frac{e(v_++v_-)}{2}\leq\pm\la\\
&\Leftrightarrow&\left(\zeta_{\nu}^2+(mc^2+\nu\frac{e(v_+-v_-)}{2})^2\right)^{\frac12}\leq\pm\left(\la-\nu\frac{e(v_++v_-)}{2}\right),
\end{eqnarray*}thus
\begin{eqnarray*}\!\int_{H_1^+(x,\xi)\leq\la}\!\!\!\!d\xi-\int_{H_1^-(x,\xi)\geq\la}\!\!\!\!d\xi \!\!\!\!&=&\!\!\!\! \frac{4\pi}{3}\!\left(\!(\la-\frac{e(v_++v_-)}{2})_{+}^2-(mc^2+\frac{e(v_+-v_-)}{2})^2\!\right)_+^{\frac{3}{2}}\\
\!\!\!\!&-&\!\!\!\!\frac{4\pi}{3}\!\left(\!(\la-\frac{e(v_++v_-)}{2})_{-}^2-(mc^2+\frac{e(v_+-v_-)}{2})^2\!\right)_+^{\frac{3}{2}},\end{eqnarray*}
 and
\[-\int_{H_0^+(x,\xi)\leq\la}\!\!d\xi+\int_{H_0^-(x,\xi)\geq\la}\!\!d\xi=\mp\frac{4\pi}{3}\left(\la^2-(mc^2)^2\right)_+^{\frac{3}{2}},\;\;\;\mbox{for}\;\pm\la\geq 0,\]
with $(x)_+=max(x,0)$ and $(x)_-=max(-x,0)$ for $x\in\R.$\hfill{$\square$}
\begin{rem}\label{remarkweyl}
Theorem \ref{thmweyl} can be extended to the pair operators $(H_1=H_0+V_1,\ H_2=H_0+V_2)$, where the potentials $V_1,V_2$ are electro-magnetic potentials 
$$V_\.(x)=e(-\al\cdot A^\.+v^\.)(x)=-\sum_{j=1}^3\al_j\cdot eA_j^\.(x)+e\left(
\begin{array}{cc}
v_+^\.(x)I_2 & 0 \\
0 & v_-^\.(x)I_2
\end{array}
\right)$$ satisfying the assumption ${\bf(A_V)}$ with $\de>0$ (or $\|V_j(x)\|\longrightarrow0$) and the potential $V=V_2-V_1$ satisfies the assumption ${\bf(A_V)}$ with $\de>3$:\\
For all $\la,\;\la_1$ noncritical energy levels for $H_1,\ H_2$ such that $\pm mc^2\not\in\; ]\la_1,\la[$ and $h\in ]0,h_0[$, 
we have the asymptotic expansion 
\begin{eqnarray}\label{weylg}
\xi(\la,h)-\xi(\la_1,h)=w(\la,\la_1)h^{-3}+O(h^{-2}).
\end{eqnarray}
Here the $O(h^{-2})$ is uniform for $\la$ (resp.$\la_1$) in a small interval $I_2$ (resp. $I_1$). The first term $w(\la,\la_1)\;\in C^\infty(I_2\times I_1)$ is given by 
$$w(\la,\la_1)=w(\la)-w(\la_1)$$ with,
\begin{eqnarray}\label{eqpremierterm1}
\;\;\;\;\;\;w(\la)=\frac{1}{3\pi^2}\int_{\R^3}\big[W_+(\la,v_+^\.,v_-^\.)-W_+(\la,v_+^\.,v_-^\.)\big]_1^2 \,dx
\end{eqnarray}
where   $\;\;W_\pm(\la,a,b)=\left(\!\left(\la-\frac{e(a+b)}{2}\right)_{\pm}^2-\left(mc^2+\frac{e(a-b)}{2}\right)^2\!\right)_+^{\frac{3}{2}}$.

In this setting, we does not have a formula like (\ref{formulae}). But it could be possible to use the approach to Bruneau-Petkov in \cite{VBVP}. For that we need more informations on the approximation of the propagator  $e^{-itH_jh^{-1}}$ by the Fourier integral operator ${\cal U}_t$.
\end{rem}

\section{ Resonances in small domains}

In this section, we 
assume that the hamiltonian $H=H_0+V$, where $V$ is an electro-magnetic potential:
\[H=-\sum_{j=1}^3\al_j(ich\p_j+eA_j)+mc^2\be+ev.\]
\subsection{Upper bound for the number of resonances in domains of width $h$}\label{majorationh} 
We adapt, for the Dirac operator, Theorem 1 of \cite{VBVP2} which rests on a representation formula for the spectral shift function (see Theorem \ref{thmxi}).
\begin{thm}\label{thmequiv}
Suppose that each $\la\in[E_0,E_1]$ is a non-critical energy level for $H$. Then for $h\in]0,h_0]$, the following assertions are equivalents:
\begin{itemize}
\item[\emph{(i)}] There exist positive constants $B,\,C,\,b,\,h_0$, such that for any $\la\in[E_0-b,E+b],$ $\;h\in]0,h_0]$ and $h/B\leq\rho\leq B$, we have
\[\# \{z\in\C:\;z\in\res(H),|z-\la|\leq\rho\}\leq C\rho h^{-3}.\] 
\item[\emph{(ii)}] There exist positive constants $B_1,C_1,\varepsilon_1,h_1,$ such that for any $\la\in[E_0-\varepsilon_1,E_1+\varepsilon_1],$ $\;h\in]0,h_1]$ and $h/B_1\leq\rho\leq B_1$, we have
\[|\xi(\la+\rho,h)-\xi(\la-\rho,h)|\leq C\rho h^{-3}.\]
\end{itemize}
\end{thm}
As a consequence of Theorem \ref{thmequiv}, we have an upper bound $O(h^{-2})$ for the number of the resonances for the semi-classical Dirac operator close to a non-critical energy level in a domain of width $h$:
\begin{prop}\label{propmagh}Assume that $V$ is the electro-magnetic potential $(\ref{eqV})$ satisfying the assumption ${\bf (A_V)}$ with $\de>3$. We
suppose also that each $\la\in[E_0,E_1]$ is a non-critical energy level for $H$. There are  positive constants $C,\,B,\,b,\,h_0$ such that for any $\la\in[E_0-b,E+b],\;h\in]0,h_0]$ and $h/B\leq\rho\leq B$, we have
\[\# \{z\in\C:\;z\in\res(H),|z-\la|\leq\rho\}\leq C\rho h^{-3}.\] 
 
\end{prop}
\Proof It follows from Theorem \ref{thmequiv} and equation (\ref{weyl}).\hfill{$\square$} 
\subsection{Breit-Wigner approximation}\label{secbreit}

In this part, we consider small domains of width $h$, and we prove a Breit-Wigner approximation for $\xi(\la,h)$ (see \cite{VPMZ1}, \cite{VPMZ2}, \cite{CGAMDR}, \cite{JBJS}, \cite{VBVP}). Let $\eta(\la,h)$ be the real analytic function defined by
\[\eta(\la,h)=\xi(\la,h)-\#\{\mu\in[E_0,\la]:\;\mu\in \si_{d}(H)\}.\]
Using Proposition \ref{propmagh} and the arguments used in \cite[Section 6]{VBVP}, we obtain a Breit-Wigner approximation for the derivative of the spectral shift function $\xi(\la,h).$
\begin{thm}\label{thmbreit}{\bf(Breit-Wigner)}
Assume that $V$ is an electro-magnetic potential ($\ref{eqV}$), for any $\la\in[E_0,E_1]$ non-critical energy level for $H$,  $0<\rho<h/B,$ $0<B_1<B$, and $h$ sufficiently small, we have 
\[\eta(\la+\rho,h)-\eta(\la-\rho,h)=\hspace{-5mm}\sum_{\begin{array}{c}
\scriptstyle w\in \res(H) \\
\scriptstyle \emph{Im}\,  w\neq0,\ |w-\la|<h/B_1
\end{array}}\hspace{-10mm}\om_{\C_-}(w,[\la-\rho,\la+\rho])+O(\rho)h^{-3},\]
where $B>0$ is  a positive constant 
and $\om_{\C_-}$ is the harmonic measure 
$$\om_{\C_-}(w,E)=-\frac{1}{\pi}\int_E\frac{\emph{Im} (w)}{|t-w|^2}dt,\;\;\;E\subset\R=\p\C_-.$$
\end{thm}
Using Theorem \ref{thmequiv} and repeating with little modifications the arguments used in  \cite[Section 6]{VBVP1}, we obtain the following corollary which implies also trace formula in small domains.
\begin{corl}Under the assumptions of Theorem \ref{thmbreit} and supposing that $[E_0,E_1]$ is a non-critical energy level for $H$, for each $E\in[E_0,E_1]$ there exist constants $C_2>C_1>0,\,h_0>0$ so that for $|\la-E|\leq C_1h,\;h\in]0,h_0],\;$we have
\begin{eqnarray}
\xi'(\la,h)=-\frac{1}{\pi}\sum_{\begin{array}{c}
\scriptstyle w\in \res(H) \\
\scriptstyle |E-w|\leq C_2h
\end{array}}\frac{\emph{Im} (w)}{|\la-w|^2}+\sum_{\begin{array}{c}
\scriptstyle  w\in \si_{d}(H) \\
\scriptstyle |E-w|\leq C_1h
\end{array}}\de_w(\la)+O(h^{-3}).
\end{eqnarray}
Here $\de_w(\cdot)$ is the Dirac 
mass at $w \in \R$.
\end{corl}
\begin{center}
\appendix{\bf Appendix A. Construction of ${\cal U}_{t}$}
\end{center}
In this appendix, we construct a parametrix at small time of the propagator of the Dirac equation in an external electro-magnetic field  $$ih\p_t \psi=H_1\psi,$$
with $H_1=H_0+V$. Here $H_0$ is the selfadjoint operator defined in (\ref{eqH_0}) and $V$ is an electro-magnetic potential (\ref{eqV}).
\begin{thm}{\bf(Approximation of the propagator)}\label{thmpropagator} 
 There exist $\de_1>0$ small enough and a Fourier integral operator  ${\cal U}_{t}={\cal U}_{t}^++{\cal U}_{t}^-$ with
$${\cal U}_{t}^\pm f(y)= \frac{1}{(2\pi h)^3}\int\!\!\int e^{i(\Phi^{\pm}(t,x,\xi)-y\cdot\xi)h^{-1}}E^{\pm}(t,x,y,\xi;h)f(y)d\xi dy,$$
defined for $|t|<\de_1$ such that:
\begin{itemize}
\item The amplitude $E^{\pm}(t,x,y,\xi;h)\in {\cal S}^0(1)$.
\item $\|{\cal U}_{t}-e^{-itH_1h^{-1}}\|_{}=O(h^{\infty}),$ uniformly for $|t|<\de_1$. 
\item The phase function $\Phi^{\pm}(t,x,\xi)-x\cdot\xi$ and their derivatives $\p_t^{\al}\p_x^{\be}\p_{\xi}^{\ga}\left(\Phi^\pm(t,x,\xi)-x\cdot\xi\right)$ are uniformly bounded for $(t,x,\xi)\in [-\de_1,\de_1]\times\R^3\times B(0,C_1) $,  $(\al,\be,\ga)\neq(0,0,0)$ and $C_1>0$ (see(\ref{eqphase})).
\end{itemize}
\end{thm}
With a different approach, such result has been obtained by Yajima \cite{KY} for a scalar electric potential ($v_+=v_-$).\\

\Proof
We consider the equivalent problem for ${\cal U}_{t}$ \begin{eqnarray}\label{equt0}\left\{\begin{array}{rl}
ih\p_t{\cal U}_{t}-H_1{\cal U}_{t}&=0,\\
{\cal U}_{0}&=I.
\end{array}\right.
\end{eqnarray}
We solve this problem using the B.K.W. method. 
We assume that the kernel of the operator ${\cal U}_{t}$ is $K_{t}$, where
\[K_{t}(x,y;h)=\frac{1}{(2\pi h)^3}\int e^{i(\Phi(t,x,\xi)-y\cdot\xi)h^{-1}}E(t,x,y,\xi;h)d\xi,\] 
with $E(t,x,y,\xi;h)=E_0(t,x,y,\xi)+hE_1(t,x,y,\xi)+\cdots$.\\

Thus, if we look for $E(t,x,y,\xi;h)$  having the asymptotic expansion above, it is enough to solve (in some fixed neighborhood of $t=0$) the sequence of equations

\begin{eqnarray}\label{eqtransport}\left\{\begin{array}{rl}
0&=\left(\p_t\Phi(t,x,\xi) +c\al\cdot\nabla_x\Phi-e\al\cdot A+mc^2\be+ev\right)E_0,\\
i(\p_t+c\al\cdot\nabla_x)E_{j}\!\!&=\!\!\left(\p_t\Phi(t,x,\xi) +c\al\cdot\nabla_x\Phi-e\al\cdot A+mc^2\be+ev\right)E_{j+1},\\
E_0(0,x,\xi)&=I_4,\;\\
E_j(0,x,\xi)&=0, \;\; \mbox{for} \;j\geq1.\;\;\;\;\;\;\;\;\;\;\;\;\;\;\;\;\;\;\;\;\;
\end{array}\right.
\end{eqnarray}
On the support of $E_0$, we deduce the eikonal equation 
\begin{eqnarray}\label{eqeikonale}\left\{\begin{array}{r}
\det\left(\p_t\Phi(t,x,\xi) +c\al\cdot\nabla_x\Phi-e\al\cdot A+mc^2\be+ev\right)=0,\\
\Phi(0,x,\xi)=x\cdot\xi.\;\;\;\;\;\;\;\;\;\;\;\;\;\;\;\;\;\;\;\;\;\;\;\;\;\;\;\;\;
\end{array}\right.
\end{eqnarray}
The system $(\ref{eqeikonale})$ is equivalent to 
\begin{eqnarray}\label{eqeiconale1}\left\{\begin{array}{r}
\p_t\Phi^\pm(t,x,\xi) +H_1^{\pm}(x,\nabla_x\Phi)=0,\;\;\;(\mbox{see}\;(\ref{eqH+-})),\\
\Phi^\pm(0,x,\xi)=x\cdot\xi.\;\;\;\;\;\;\;\;\;\;\;\;\;\;\;\;\;\;\;\;\;\;\;\;\;\;\;\;\;
\end{array}\right.
\end{eqnarray}
The last system can be solved by Hamilton-Jacobi method (see \cite{VA}) and all derivatives \begin{eqnarray}\label{eqphase}
\p_t^{\al}\p_x^{\be}\p_{\xi}^{\ga}\left(\Phi^\pm(t,x,\xi)-x\cdot\xi\right)
 \end{eqnarray}
 are uniformly bounded for $(t,x,\xi)\in [-\de_1,\de_1]\times\R^3\times B(0,C_1) $ and $(\al,\be,\ga)\neq(0,0,0)$.\\
Using the Taylor formula in a neighborhood of $t=0$, the two solutions of $(\ref{eqeiconale1})$ satisfy:
 \[\Phi^{\pm}(t,x,\xi)= x\cdot\xi-tH_1^{\pm}(x,\xi)+O(t^2).\]
 Then 
 ${\cal U}_{t}={\cal U}_{t}^++{\cal U}_{t}^-$,
and the kernel of the operator ${\cal U}_{t} $ is $K_{t}=K_{t}^++K_{t}^-$, with
\[K_{t}^{\pm}(x,y;h)=\frac{1}{(2\pi h)^3}\int e^{i(\Phi^{\pm}(t,x,\xi)-y\cdot\xi)h^{-1}}E^{\pm}(t,x,y,\xi;h)d\xi.\] 
We look for the amplitude $E^{\pm}(t,x,y,\xi;h)$ having an asymptotic expansions according to the power of $h$:$$E_0^{\pm}(t,x,y,\xi)+hE_1^{\pm}(t,x,y,\xi)+\cdots.$$ Consequently, the coefficients $E_j^{\pm}(t,x,y,\xi)$ are the solutions of the transport equations

\begin{eqnarray}\label{eqtransport1}\left\{\begin{array}{rl}
0=&\left(\p_t\Phi^{\pm}+c\al\cdot\nabla_x\Phi^{\pm}-e\al\cdot A+mc^2\be+ev\right)E_{0}^{\pm},\\
i(\p_t+c\al\cdot\nabla_x)E_{j}^{\pm}=&\left(\p_t\Phi^{\pm}+c\al\cdot\nabla_x\Phi^{\pm}-e\al\cdot A+mc^2\be+ev\right)E_{j+1}^{\pm},\\
E_j^+(0,x,\xi)+E_j^-(0,x,\xi)&=0 \;\;\; \mbox{for} \;j\geq1,\;\;\;\;\;\;\;\;\;\;\;\;\;\;\;\;\;\;\;\\
E_0^{\pm}(0,x,\xi)&=\Pi_1^{\pm}(x,\xi),\;
\end{array}\right.
\end{eqnarray}
with $\Pi_1^{\pm}(x,\xi)$ defined by $(\ref{eqP+-}).$\\ \\
{\bf Resolution of (\ref{eqtransport1}).}\\
Let us denote by $L=\p_t+c\al\cdot\nabla_x$, with $\al\cdot\nabla_x=\sum_{j=1}^3\al_j\p_{x_j}.$
The matrix $${\MM}^{\pm}=\p_t\Phi^{\pm}+c\al\cdot\nabla_x\Phi^{\pm}-e\al\cdot A+mc^2\be+ev,$$ is Hermitian and has two real eigenvalues linearly independent with multiplicity $2$. \\

 First, we multiply the system $(\ref{eqtransport1})$ by the column-vector $N_1=(1,0,0,0)^{\dagger}$, the subscript $\dagger$ designates the complex conjugate of the transpose. We denote \begin{eqnarray}\label{supposition}
 E_{j,1}^{\pm}=E_j^{\pm}N_1\; \mbox{for} \;j=1,2,\cdots,\;\;\; E_{0,1}^{\pm}(0,x,\xi)=\Pi_1^{\pm}(x,\xi)N_1.
 \end{eqnarray}
  
 Since $\det({\MM}^{\pm})=0$, there exist $l_k^{\pm}$ and $r_k^{\pm}$, left and right eigenvectors of the  matrix ${\MM}^{\pm}$, corresponding to the eigenvalue zero, such that \begin{eqnarray}\label{eqnoyau}{\MM}^{\pm}r_k^{\pm}=0,\;\;\;\;\;l_k^{\pm}{\MM}^{\pm}=0,\;\;\;\;\;\;
l_k^{\pm}=(r_k^{\pm})^{\dagger},\;\;\;\;\;\;k=1,2,
\end{eqnarray} 
(here $r_k^{\pm}$ is a column-vector and $l_k^{\pm}$ is a row-vector). We choose
\[r_1^{+}=\left(\begin{array}{c}
u^{+}\\
0\\
v^{+}\\
w_+^{+}\end{array}\right),\;\;\;r_2^{+}=\left(\begin{array}{c}
0\\
u^{+}\\
w_-^{+}\\
-v^{+}\end{array}\right),\;\;r_1^{-}=\left(\begin{array}{c}
w_+^-\\
v^-\\
0\\
u^-\end{array}\right),\;\;r_2^{-}=\left(\begin{array}{c}
-v^-\\
w_-^-\\
u^-\\
0\end{array}\right),\]
\begin{eqnarray}\label{eqvp}
l_{\nu}^{\pm}r_k^{\pm}=(\mp2p_5^{\pm}u^{\pm})\de_{\nu k},\;\;\;\;\;\;\nu,k=1,2.
\end{eqnarray}
Here $u^{\pm},v^{\pm}$ and $w_{\pm}^{\pm}$ are defined by 
\[u^{\pm}=p_4^{\pm}\mp p_5^{\pm},\;\;\;v^{\pm}=p_3^{\pm},\;\;\;w_+^{\pm}=\pm p_1^{\pm}+ip_2^{\pm},\;\;\;w_-^{\pm}=\pm p_1^{\pm}-ip_2^{\pm},\]
where $ p_4^{\pm}=mc^2+\frac{e(v_+-v_-)}{2},\;\;\;\;p_5^{\pm}=\p_t\Phi^{\pm}+\frac{e(v_++v_-)}{2},\;\;\;p_j^{\pm}=c\p_{x_j}\Phi^{\pm}-eA_j,$ for $j=1,2,3.$

It is easy to see that the vector-valued functions $r_k^{\pm}(t,x,\xi)$ and $l_k^{\pm}(t,x,\xi)$ can be chosen to be smooth in $t$ and $x$ and different from $0$ everywhere. All the derivatives of $r_k^{\pm},l_k^{\pm},$  $k=1,2,$ are uniformly bounded for $(t,x,\xi)\in [-\de_1,\de_1]\times\R^3\times B(0,C_1)$. Then it follows from the first equation in ($\ref{eqtransport1}$) that $$E_{0,1}^{\pm}=\si_{0,1}^{\pm}(t,x,\xi)r_1^{\pm}(t,x,\xi)+\si_{0,2}^{\pm}(t,x,\xi)r_2^{\pm}(t,x,\xi),$$ where $\si_{0,1}^{\pm}, \si_{0,2}^{\pm}$, are scalar-valued functions. If we multiply the second equation in $(\ref{eqtransport1})$ for $j=0$ on the left by $l_k^{\pm}$ for $k=1,2$, 
we deduce the following differential equations  for $\si_{0,k}^{\pm}$:$$\left\{\begin{array}{rl}
l_1^{\pm}L(\si_{0,1}^{\pm}r_1^{\pm})+l_1^{\pm}L(\si_{0,2}^{\pm}r_2^{\pm})&=0,\\
l_2^{\pm}L(\si_{0,1}^{\pm}r_1^{\pm})+l_2^{\pm}L(\si_{0,2}^{\pm}r_2^{\pm})&=0.
\end{array}
\right.
 $$
We conclude
\begin{eqnarray}\label{eqs}\left\{\begin{array}{rl}
l_1^{\pm}r_1^{\pm}\p_t(\si_{0,1}^{\pm})+c\sum_{j=1}^3l_1^{\pm}\al_jr_1^{\pm}\p_{x_j}(\si_{0,1}^{\pm})+c\sum_{j=1}^3l_1^{\pm}\al_jr_2^{\pm}\p_{x_j}(\si_{0,2}^{\pm})\\ \\ +l_1^{\pm}L(r_1^{\pm})\si_{0,1}^{\pm}+l_1^{\pm}L(r_2^{\pm})\si_{0,2}^{\pm}=0,\\ \\ 
l_2^{\pm}r_2^{\pm}\p_t(\si_{0,2}^{\pm})+c\sum_{j=1}^3l_2^{\pm}\al_jr_2^{\pm}\p_{x_j}(\si_{0,2}^{\pm})+c\sum_{j=1}^3l_2^{\pm}\al_jr_1^{\pm}\p_{x_j}(\si_{0,1}^{\pm})\\ \\ +l_2^{\pm}L(r_1^{\pm})\si_{0,1}^{\pm}+l_2^{\pm}L(r_2^{\pm})\si_{0,2}^{\pm}=0.
\end{array}
\right.
\end{eqnarray}
We now use Lemma $\ref{lemsi}$ (see below) in the system ($\ref{eqs}$). Since $p_5^{\pm}\neq0$, $u^{\pm}=p_4^{\pm}\mp p_5^{\pm}\neq0$ then, after multiplying ($\ref{eqs}$) by $(\mp 2p_5^{\pm}u^{\pm})^{-1}$, the system ($\ref{eqs}$) can be written as 
\begin{eqnarray}
D^{\pm}\si_0^{\pm}=M^{\pm}\si_0^{\pm}:=(\mp 2p_5^{\pm}u^{\pm})^{-1}\left(\begin{array}{cc}
l_1^{\pm}L(r_1^{\pm})&l_1^{\pm}L(r_2^{\pm})\\
l_2^{\pm}L(r_1^{\pm})&l_2^{\pm}L(r_2^{\pm})
\end{array}
\right)\left(\begin{array}{c}
\si_{0,1}^{\pm}\\
\si_{0,2}^{\pm}
\end{array}
\right),
                                                \end{eqnarray}
 with $D^{\pm}=\p_t+a^{\pm}\cdot\nabla_x=\p_t+\sum_{j=1}^3a_j^{\pm}(t,x)\p_{x_j},$ and\\ $$\;\;a^{\pm}=c(\mp2p_5^{\pm}u^{\pm})^{-1}(l_1^{\pm}\al_1r_1^{\pm},\;\;l_1^{\pm}\al_2r_1^{\pm},\;\;l_1^{\pm}\al_3r_1^{\pm}).$$
Thus the function $\si_{0,k}^{\pm}$ can be found if its value is known for $t=0$, and it is as smooth as $\si_{0,k}^{\pm}(0,x,\xi)$ (for more details, see a method 
for solving a similar equation in \cite{SRJK}). The equality 
$$E_{0,1}^{\pm}(0,x,\xi)=\si_{0,1}^{\pm}(0,x,\xi)r_1^{\pm}(0,x,\xi)+\si_{0,2}^{\pm}(0,x,\xi)r_2^{\pm}(0,x,\xi)=\Pi_1^{\pm}(x,\xi)N_1,$$ 
gives the value of $\si_0^{\pm}$ at $t=0$.\\
Since, all the derivatives of $\si_{0,k}^{\pm},\;r_k^{\pm},$ for $\;k=1,2,$ are uniformly bounded, then all the derivatives $(\p_t^{\al}\p_x^{\be}\p_{\xi}^{\ga}E_{0,1}^{\pm})$ are uniformly bounded for $(\al,\be,\ga)\in \N\times\N^3\times\N^3.$\\ \\
It follows from the second equation in $(\ref{eqtransport1})$ for $j=0$, that
\[iLE_{0,1}^{\pm}={\MM}^{\pm}E_{1,1}^{\pm},\]
i.e., $E_{1,1}^{\pm}=\si_{1,1}^{\pm}r_1^{\pm}+\si_{1,2}^{\pm}r_2^{\pm}+h_1^{\pm},$ where $\si_{1,k}^{\pm}$ is a scalar-valued function for $k=1,2,$ and $h_1^{\pm}$ is expressed in term of $LE_{0,1}^{\pm}$. To find $\si_{1,k}^{\pm}$ it is 
sufficient to multiply the second equation in $(\ref{eqtransport1})$ for $j=1$ on the left by $l_k^{\pm}$ for $k=1,2$. Then
$$\left\{\begin{array}{rl}
l_1^{\pm}L(\si_{1,1}^{\pm}r_1^{\pm})+l_1^{\pm}L(\si_{1,2}^{\pm}r_2^{\pm})+l_1^{\pm}L(h_1^{\pm})&=0,\\ \\
l_2^{\pm}L(\si_{1,1}^{\pm}r_1^{\pm})+l_2^{\pm}L(\si_{1,2}^{\pm}r_2^{\pm})+l_2^{\pm}L(h_1^{\pm})&=0.
\end{array}
\right.
 $$

From this equation, $\si_{1,k}^{\pm}$ can be found if the function $\si_{1,k}(0,x)$ is known. With the same process, for all $j=1,2,\cdots$, we obtain
$$\left\{\begin{array}{rl}
\si_{j,1}^{\pm}r_1^{\pm}+\si_{j,2}^{\pm}r_2^{\pm}+h_j^{\pm}&=E_{j,1}^{\pm},\\
l_1^{\pm}L(\si_{j,1}^{\pm}r_1^{\pm})+l_1^{\pm}L(\si_{j,2}^{\pm}r_2^{\pm})+l_1^{\pm}L(h_j^{\pm})&=0,\\
l_2^{\pm}L(\si_{j,1}^{\pm}r_1^{\pm})+l_2^{\pm}L(\si_{j,2}^{\pm}r_2^{\pm})+l_2^{\pm}L(h_j^{\pm})&=0.
\end{array}
\right.
 $$
For $t=0$, $\;\;\;j=1,2,\cdots$, we have  $$\si_{0,1}^{\pm}r_1^{\pm}+\si_{0,2}^{\pm}r_2^{\pm}=\Pi_1^{\pm}N_1,\;\;\;\si_{j,1}^+r_1^++\si_{j,1}^-r_1^-+\si_{j,2}^+r_2^++\si_{j,2}^-r_2^-=-(h_j^++h_j^-),$$
and the quantity $h_j^{\pm}$ is determined if $E_{0,1}^{\pm},E_{1,1}^{\pm},\cdots,E_{j-1,1}^{\pm}$, are known. 
Solving the differential equation for $\si_j^{\pm}=\left(\begin{array}{c}
\si_{j,1}^{\pm}\\
\si_{j,2}^{\pm}\end{array}\right)$, we find these functions for all sufficiently small $t$.\\ \\
Repeating this block of calculus, multiplying by $N_2=(0,1,0,0)^{\dagger},\;N_3=(0,0,1,0)^{\dagger}$ and $N_4=(0,0,0,1)^{\dagger}$ instead of $N_1$ in $(\ref{supposition})$, we find $E_{j,2}^{\pm}\!=E_j^{\pm}N_2,\;E_{j,3}^{\pm}\!=E_j^{\pm}N_3$ and $E_{j,4}^{\pm}\!=E_j^{\pm}N_4$.\\
Consequently, we have:
 \begin{prop}
 There exists  a family of matrices $$E_j^{\pm}=(E_{j,1}^{\pm},\;E_{j,2}^{\pm},\;E_{j,3}^{\pm},\;E_{j,4}^{\pm}),\;\;\;\;\mbox{for}\;j\geq0,$$ solution of $(\ref{eqtransport1})$. Moreover, for all $j\geq0$, $E_j^{\pm}\in C^{\infty}$ and all derivatives $(\p_t^{\al}\p_x^{\be}\p_{\xi}^{\ga}E_j^{\pm})$ are uniformly bounded for all $(t,x,\xi)\in[-\de_1,\de_1]\times\R^3\times B(0,C_1)$ and $(\al,\be,\ga)\in \N\times\N^3\times\N^3.$
 \end{prop}
 Consequently, Borel process provides a symbol $E^{\pm}(t,x,y,\xi;h)\in {\cal S}^0(1)$ with compact support on $\xi$ and $(x-y)$ with $E_0^{\pm}(t,x,y,\xi)+hE_1^{\pm}(t,x,y,\xi)+\cdots$ its asymptotic expansion.
 \\\\
{\bf Desired estimate.}\\
Next, we remark that for all $N\in\N$:
\begin{eqnarray}(ih\p_t-H_1)\left(e^{i(\Phi^{\pm}(t,x,\xi)-y\cdot\xi)h^{-1}}\sum_{j=0}^{N}h^jE_j^{\pm}\right)=&\!\!\!\!\!\!\!\!\!\!\!\!\!\!\!\!\!\!\!\!\!\!\!\!\!\!\!\!\!\!\!e^{i(\Phi^{\pm}(t,x,\xi)-y\cdot\xi)h^{-1}}\nonumber\\
\times& \sum_{j=0}^{N}(ihL(E_j^{\pm})+{\MM}^{\pm}E_{j}^{\pm})h^j\nonumber\\
=&\!\!\!\!\!\!\!\!\!\!\!\!\!\!\!\!\!\!\!\!\!\!\!\!\!\!\!\!\!\!\!\!\!\!P_N(t,x,\xi;h)h^N,
\end{eqnarray}
and all derivatives $D_{x,\xi}^\al P_N(t,x,\xi;h)$ are bounded as $h\to0$ for all $\al$.
Then for all $N\in \N$,  
\begin{eqnarray}\left\{\begin{array}{rl}
ih\p_t{\cal U}_{t}-H_1{\cal U}_{t}&=O(h^N),\\
{\cal U}_{0}&=I+O(h^N),
\end{array}\right.
\end{eqnarray}
thus\begin{eqnarray}\left\{\begin{array}{rl}
\frac{d}{dt}(e^{+itH_1h^{-1}}{\cal U}_{t})&=O(h^N),\\
{\cal U}_{0}&=I+O(h^N),
\end{array}\right.
\end{eqnarray}
where $O(h^N)$ is uniform on $t$ and corresponds to the norm in ${\cal L}(L^2)$. Then 
we get:
\begin{eqnarray}
\label{eqnorm}\|{\cal U}_{t}-e^{-itH_1h^{-1}}\|=O(h^{\infty}).
\end{eqnarray}

\hfill{$\square$}
\begin{lem}\label{lemsi}Under the notations used above, we have 
\begin{eqnarray}\label{eqr} l_1^{\pm}\al_jr_1^{\pm}=l_2^{\pm}\al_jr_2^{\pm},\;\;\;\;\;\;l_1^{\pm}\al_jr_2^{\pm}=l_2^{\pm}\al_jr_1^{\pm}=0, \;\;j=1,2,3.
\end{eqnarray}
\end{lem}
\Proof
As Rubinow and Keller in \cite{SRJK}  let us work in a general situation.\\ 
We consider the $n$ Hermitian matrices $M_{\mu}$ and $n$ real scalars $p_{\mu}$, $\mu=1,\cdots,n.$ 
Let $G$ be the Hermitian matrix defined by$$G=\sum_{\mu=1}^np_{\mu}M_{\mu}.$$
Let $\la$ be a multiple eigenvalue of $G$ and $B_1,\cdots,B_q,$ a set of corresponding orthormal eigenvectors which are differentiable functions of $p_{\mu}$. Then \begin{eqnarray}
B_j^{\dagger}B_k&=&\de_{jk},\\
 GB_k&=&\la B_k\label{eqla}.   
\end{eqnarray}
The subscript $\dagger$ designates the complex conjugate of the transpose.\\
If $\la(p_1,\cdots,p_{\mu})$ is differentiable, we differentiate $(\ref{eqla})$ with respect to $p_{\mu}$ and obtain
\begin{eqnarray}\label{eqdiff}
M_{\mu}B_k+G\frac{\p B_k}{\p p_{\mu}}=\frac{\p \la}{\p p_{\mu}}B_k+\la\frac{\p B_k}{\p p_{\mu}}.
\end{eqnarray}
 The multiplication of ($\ref{eqdiff}$) on the left by $B_j^{\dagger}$, the use of ($\ref{eqla}$), and the fact that $G$ is Hermitian yield
\begin{eqnarray}\label{eqbj}
B_j^{\dagger}M_{\mu}B_k=\frac{\p \la}{\p p_{\mu}}\de_{jk}.
\end{eqnarray}
In order to treat our case, we take$$G= {\cal M}^{\pm}=\sum_{\mu=1}^5p_{\mu}^{\pm}M_{\mu},$$
where $M_{j}=\al_{j}$ for $j=1,2,3$, $M_4=\be$ and $M_5=I_4$ are Hermitian matrices ($\al_j,\be$ are the Dirac matrices) and $p_{\mu}^{\pm}$ are five real scalars.\\ We also take $\la^{\pm}=p_5^{\pm}\pm\sqrt{(p_1^{\pm})^2+(p_2^{\pm})^2+(p_3^{\pm})^2+(p_4^{\pm})^2}$ and   $F^{\pm}$ the point which have the coordinate  $p_{\mu}^{\pm}$: $p_j^{\pm}=c\p_{x_j}\Phi^{\pm}-eA_j$ for $j=1,2,3$, $p_4^{\pm}=mc^2+\frac{e(v_+-v_-)}{2},\;p_5^{\pm}=\p_t\Phi^{\pm}+\frac{e(v_++v_-)}{2}$. \\ 

When $\Phi^\pm$ satisfies ($\ref{eqeikonale}$) and ($\ref{eqeiconale1}$), $r_1^{\pm},r_2^{\pm}$ are two orthogonal eigenvectors of ${\cal M}^{\pm}$ corresponding to the eigenvalue $\la^{\pm}=\la^{\pm}(F^{\pm})=0$. Since $|e(v_+-v_-)|<2mc^2$ (see ($\ref{eqv+v-}$)), $\la^{\pm}$ is differentiable near the point $F^{\pm}$. Now, we apply ($\ref{eqbj}$) with $B_j^{\dagger}=l_j^{\pm}$ and $B_k=r_k^{\pm}$. After the normalization of $r_k^{\pm},l_j^{\pm}$ we obtain
\[l_j^{\pm}M_{\mu}r_k^{\pm}=\frac{\p \la^{\pm}\left(p_1^{\pm},\cdots,p_5^{\pm}\right)}{\p p_{\mu}^{\pm}}\Big|_{F^{\pm}}(\mp2p_5^{\pm}u^{\pm})\de_{jk},\]
 and we get the lemma.\hfill{$\square$}  \\\\
\emph{Acknowledgments}. The author is grateful to V. Bruneau and J.-F. Bony for many helpful discussions. We also thank the French ANR (Grant no. JC0546063) for the partial financial support.

\end{document}